\documentclass[11pt]{article}
\usepackage{amssymb}
\usepackage{amsthm}
\usepackage{amsmath,amsfonts,amssymb}

\usepackage[mathscr]{eucal}
\usepackage{exscale}
\usepackage{natbib}
\usepackage{bm}
\usepackage{eqlist}
\usepackage[dvipsnames]{color}
\usepackage[Lenny]{fncychap}
\usepackage{multirow}
\usepackage{graphicx}
\usepackage{algpseudocode}
\usepackage{algorithm}
\usepackage{hyperref}
\usepackage{framed}

\usepackage{tabu}

\newtheorem{theorem}{Theorem}[section]

\newtheorem{remark}[theorem]{Remark}

\font\bigbf=cmbx10 scaled \magstep3

\begin{document}

\title{\bigbf Nonlinear Decision Rule Approach for Real-Time Traffic Signal Control for Congestion and Emission Reductions}

\author{Junwoo Song$^{a}\thanks{e-mail: j.song15@imperial.ac.uk}$
 \quad Simon Hu$^{a,b}\thanks{e-mail: j.s.hu05@imperial.ac.uk}$
  \quad Ke Han$^{a}\thanks{Corresponding author, e-mail: k.han@imperial.ac.uk}$
 \quad Chaozhe Jiang$^{c}\thanks{e-mail: jiangchaozhe@swjtu.edu.cn}$ 
\\\\
 $^{a}$\textit{\small Department of Civil and Environmental Engineering,}\\
\textit{\small Imperial College London, London SW7 2BU, UK}
\\
$^{b}$\textit{\small School of Civil Engineering,}
\\
\textit{\small ZJU-UIUC Institute, Zhejiang University, Haining, China}
\\
$^{c}$\textit{\small School of Transportation and Logistics,}
\\
\textit{\small Southwest Jiaotong University, Chengdu, China}
}

\maketitle

\begin{abstract}
We propose a real-time signal control framework based on a nonlinear decision rule (NDR), which defines a nonlinear mapping between network states and signal control parameters to actual signal controls based on prevailing traffic conditions, and such a mapping is optimized via off-line simulation. The NDR is instantiated with two neural networks: feedforward neural network (FFNN) and recurrent neural network (RNN), which have different ways of processing traffic information in the near past, and are compared in terms of their performances. The NDR is implemented within a microscopic traffic simulation (S-Paramics) for a real-world network in West Glasgow, where the off-line training of the NDR amounts to a simulation-based optimization aiming to reduce delay, CO$_2$ and black carbon emissions. The emission calculations are based on the high-fidelity vehicle dynamics generated by the simulation, and the AIRE instantaneous emission model. Extensive tests are performed to assess the NDR framework, not only in terms of its effectiveness in reducing the aforementioned objectives, but also in relation to local vs. global benefits, trade-off between delay and emissions, impact of sensor locations, and different levels of network saturation. The results suggest that the NDR is an effective, flexible and robust way of alleviating congestion and reducing traffic emissions.
\end{abstract}

\noindent {\it Keywords}: real-time signal control; nonlinear decision rule; congestion; emissions; neural networks

\section{Introduction}\label{secIntro}

Urban traffic signal controls play an essential role for traffic management to reduce congestion and alleviate adverse environmental impacts. Different traffic signal control strategies have been developed and deployed in large scale in the past several decades \citep{Sunkari}, ranging from traditionally pre-timed signal control systems based on historical traffic information to fully responsive systems that frequently update signal control parameters and/or phasing schemes according to real-time traffic conditions. Some typical examples of the latter include SCOOT \citep{SCOOT}, SCAT \citep{SCAT}, OPAC \citep{OPAC}, PRODYN \citep{PRODYN}, TRANSYT and RHODES.

 In the real world, traffic flows may vary significantly at road intersections even in the same time period of the day and day of the week. As a result, the capability to handle uncertain flow patterns on a network level is crucial in the design of adaptive signal controls \citep{Yin2008, PS2015, LHGFY2015, HHD2014, FHKZ2015, CAS2016, SBXR2013}. The objectives of adaptive signal control strategies include minimization of (weighted) vehicle/pedestrian delay \citep{HHD2014, SBW2006, ZYL2010}, minimization of passenger delay \citep{CS2011, CAS2016}, minimization of number of stops \citep{LMH2000}, maximization of total throughput \citep{CS2004, HGPFY2014}. 
 
This paper focuses on real-time adaptive signal control on realistic traffic networks, while taking into account exhaust emissions including total carbon and Black Carbon (BC). Total carbon is closely related to the emission of CO2, the primary greenhouse gas contributing to the climate change. BC is produced through incomplete combustion of carbonaceous materials, and causes serious health concerns such as respiratory problems, heart attacks and lung caners \citep{USEPA2012, JGLSCHFBK2013}. It is known that the total carbon emissions are highly dependent on the engine load and vehicle speed, while emissions of BC and NOx are more sensitive to vehicle dynamics (such as acceleration and idle) and vehicle technology \citep{ZBD2011}. Therefore, to accurately account for these different emission mechanisms in a dynamic and uncertain control environment poses a significant challenge.

The accurate modeling of different species of exhaust emissions requires high-fidelity and high-resolution traffic model and data, which provide detailed and critical information on vehicle speed, acceleration, deceleration, fleet composition, and emission factors. However, the computational burdens associated with these models typically render real-time and large-scale application of signal control and optimization infeasible. On the other hand, the decentralization of controls, in which the signal control parameters are determined at individual intersections, offers viable solutions but do not guarantee global optimality due to the lack of coordination.

In seek of a global optimal, real-time signal control strategy with multiple objectives including vehicle emissions with high fidelity and resolution, this paper proposes a novel nonlinear decision rule (NDR) approach based on feedforward neural network and recurrent neural network. The key novelty is that all the expensive computations are performed in an off-line environment through simulation-based optimization based on traffic microsimulation (S-Paramics) and high-fidelity emission modeling using AIRE and COPERT IV models \citep{MHHNVTBL2016}. The aim of the off-line optimization is to train the NDR such that its on-line (i.e. real-time) operation can be continuously improved. In addition, the on-line operation of the NDR is computationally efficient as all the optimizations are performed off line. As we shall see later, some other advantages of this framework include:
\begin{itemize}
\item {\bf flexible input structure}: The system can accommodate a wide range of data types, spatial coverage and temporal resolution. This is a desirable feature for real-time signal control as most existing studies assume full knowledge of traffic states at all key intersections and their approaches, which is often not the case in real-world networks. As shown in the case study, 
\item {\bf flexible scope and resolution of controls}: Different signal parameters (cycle time, green split, offset) at one or several intersections can be controlled simultaneously in real time;
\item {\bf user defined objectives and priorities}: As the training of the NDR is based on simulation, the proposed framework can include various traffic and environmental performance indicators; and
\item {\bf explicit incorporation of uncertainties}: Demand variations and uncertainties inherent in traffic dynamics can be accounted for during the training of the NDR, so that the resulting real-time controls are robust agains traffic uncertainties.
\end{itemize}

 While there has been numerous studies applying artificial intelligence models (such as neural networks, reinforcement learning) to real-time traffic signal controls \citep{ALUK2010, SKD2015, CHM2017}, few have considered realistic traffic dynamics that give rise to accurate estimation of emissions on a network scale. This paper employs a microscopic traffic simulation model (S-Paramics), which has been thoroughly calibrated for a real-world test network \citep{MHHNVTBL2016}, and the AIRE instantaneous emission model to accurately calculate emissions in view of different vehicle dynamics and fleet composition.

 Furthermore, the nonlinear relationships between traffic and environmental performance indicators have not been explored fully in the literature. This paper investigates the potential trade-off traffic and environmental objectives both globally (network level) and locally (junction level), as well as for different degrees of network saturation. Moreover, the impact of sensor locations on the performance of the signal controls are assessed. Our findings provide valuable insights into the management of dynamic and complex traffic networks with environmental considerations.

The proposed framework is tested for a real-world network located in Glasgow, Scotland, as part of the CARBOTRAF project \citep{MHHNVTBL2016}. Simulation-based validation of the signal controls in a real-time environment indicates a reduction of network-wide delay by up to 68\%, total carbon and black carbon emissions by 3\% and 2\%, respectively, and 1\% increase of network throughput. It is found that most of the emission reductions are concentrated at signal intersections, where local improvements can be up to 30\%. In addition, it is shown that CO$_2$ reductions, which took place primarily around traffic intersections, are correlated to delay reductions, while such correlation is weak for black carbon due to other factors like stop-and-go cycles and vehicle composition that contribute to BC emissions. Finally, the proposed NDR framework is tested with a different spatial configuration of traffic sensors, showing its robustness against sensor locations, which is a desirable property for real-world implementations.

 The rest of this paper is organized as follows. Section \ref{secLR} offers an overview of real-time signal controls in the literature. Section \ref{secFramework} outlines the general model for the NDR approach as well as implementation details of its components. Section \ref{secCaseStudy} details a case study of a real-world network and demonstrates the effectiveness of the proposed control strategies. Finally, Section \ref{secConclusion} offers some concluding remarks.

 \section{Related Work}\label{secLR}


In real-world traffic networks, traffic demands may vary significantly even in the same time period of the day and day of the week. As a result, the capability to handle stochastic flow patterns while maintaining a sound performance on a network level is crucial for the design of effective signal controls. Numerous studies are dedicated to designing adaptive or robust signal control algorithms or systems. \cite{Yin2008} develops a pre-timed signal control model by aiming to minimize the average delay and maintain sound performance against the worst-case scenario. On the network-wide level, \cite{LHGFY2015} propose a linear decision rule approach for real-time signal control. The linear decision rule relies on closed-form transformation from the state space to the control space, which is feasible in a real-time decision environment. Such a transformation can be trained via an off-line procedure, which is formulated as a distributionally robust optimization problem. \cite{HLGFY2016} propose a MILP approach to optimize signal timings that reduce network congestion as well as vehicle emissions. The MILP is developed using a robust optimization approach based on a macroscopic approximation of the relationship between link dynamics and emission rates.  \cite{ZYL2010} consider daily variations of the traffic demand in the optimization of pre-timed signal controls, by using a stochastic programming model that is informed by a range of demand scenarios and their corresponding probabilities of occurrence. \cite{URP2010} proposes a robust system optimal signal control model with an embedded cell transmission model, to account for uncertainty of future transportation demand.

From an optimization point of view, a well-defined function is required to relate the signal parameters to specific objective being optimized. As mentioned earlier, specific objectives in the literature include the minimization of (weighted) vehicle/pedestrian delay \citep{HHD2014, SBW2006, ZYL2010}, minimization of passenger delay \citep{CS2011, CAS2016}, minimization of number of stops \citep{LMH2000}, maximization of total throughput \citep{CS2004, HGPFY2014}. Furthermore, there are also numerous studies that incorporate environmental objectives such as emission and fuel consumption. \cite{HLGFY2016} propose a signal optimization method that takes advantage of a macroscopic relationship between link occupancy and vehicle emission rate. Their study is based on the Lighthill-Whitham-Richards kinematic wave model, from which vehicle-derived emissions are calculated. Through robust optimization, the authors are able to reformulate signal optimization problems with emission constraints/objectives as a mixed integer linear program. \cite{JHHT2014} also have developed a method to optimize transit signal priority scheme by alleviating impact on exhaust emission and reducing traffic vehicle delay. However, it finds, in many cases, traffic and emission objectives are not aligned very well with each other, especially when traffic network is complicated and traffic dynamics are nonlinear. Thus, in order to keep trade-off between both objectives, developing a bi-objective optimization model for traffic signal setting has gained popularity. \cite{SSSO2015} propose a novel integration method in order to solve multi-objective traffic signal optimization. The method can keep a balance between mobility, safety, and exhaust emission by communicating between models in the integration method. \cite{CBMW2012} mentioned that vehicle emissions are affected by a variety of factors; vehicle type, vehicle operation time and condition (idle speed, acceleration, and deceleration). So, the instantaneous vehicle emission model based on detailed vehicle dynamics is more appropriate for this kind of study. Most of these aforementioned signal optimization strategies rely on either simplified vehicle dynamics (such as the kinematic wave model) or fleet composition (e.g. single commodity). It is widely known that an accurate depiction of traffic emissions requires extensive knowledge of the detailed vehicle movements, vehicle types, as well as relevant emission factors \citep{MHHNVTBL2016}. However, such information is very difficult to obtain especially on a network-wide scale during a real-time operational environment, and most signal optimization algorithms tend to resort to heuristics. In addition, the potential trade-off between traffic performance and environmental impact has not been properly understood in an on-line decision-making context.

The environmental impact of traffic signal control strategies has been investigated and accounted for in a number of recent studies. \cite{LDXH2013} consider vehicle mean speed and the number of vehicle stops to simultaneously reduce vehicle delay and traffic emissions for urban traffic networks by applying model predictive control (MPC). Similarly, \cite{JPPD2017} use the MPC  with a gradient-based control optimization approach to smooth vehicle flows, in order to reduce traffic congestion and emissions simultaneously. \cite{CH2016} develop a traffic emission control model considering signal timing and emission pricing. The model, based on particle swarm optimization, is able to optimize intersection traffic and link-based emissions. \cite{ZYC2013} use the cell transmission model to describe traffic dynamics and vehicle emissions, and devise a signal optimization scheme that takes into account pollutant dispersion affected by weather conditions. \cite{ZC2014} develop a multi-objective optimization method based on microscopic traffic simulation at a single intersection, A modal emission and fuel consumption model is used in conjunction with the genetic algorithm to minimize vehicle delay, exhaust emission and fuel consumption at the same time. \cite{OSO2015} propose a meta-model, simulation-based approach to optimize fixed timing for dynamic traffic networks by incorporating dynamic traffic assignment models. The response surface methodology is shown to significantly reduce the computational burden typically associated with microscopic traffic and emission models.

Recent advancement in artificial intelligence, both in theory and computational architecture, has led to the emergence of a number of {\it machine learning} (ML) based approaches for traffic signal controls, such as {\it neural networks} (NNs) and {\it reinforcement learning} (RL). In particular, multi-agent approach using NN has been applied to minimize average vehicular delay time and average stoppage time \citep{SCC2006}, improve the reactivity of traffic control and capacity of traffic network \citep{CHM2017}, alleviate traffic congestion \citep{SBH2013}, optimize driving policies \citep{Wiering2000}, and improve traffic control decision-making \citep{HS2001, SKD2015}. On the other hand, RL is used to develop multi-agent traffic control architecture to optimize phase timing \citep{BGS2010}, reduce queue length and the number of stops \citep{LLW2016}, and minimize the average delay and congestion at intersections \citep{ALUK2010, CWH2009}.

\section{NDR Based Control Framework}\label{secFramework}

The nonlinear decision rule (NDR) approach for real-time signal control problems is detailed in this section. In presenting the model we first employ a generic representation without relying on any specific network configuration or control preferences, which highlights the flexibility and robustness of the proposed method. This is done in Section \ref{subsecgenericmodel}. Implementation details of the model pertaining to the case study of this paper will be presented in Section \ref{subsecimplementation}. Finally, the off-line training of the NDR based on simulation-based optimization will be detailed in Section \ref{subsecofflinetraining}.

\subsection{The generic model}\label{subsecgenericmodel}

The dynamics of the traffic network of interest may be perceived by a state vector $\bm q=\bm q(t)$ that changes with time. For example, the vector $\bm q$ may be used to express traffic quantities such as flow, density, speed, and travel time, which may be measured with different types of sensors (e.g. loop detectors, GPS, and cameras). In addition, we allow $\bm q$ to encapsulate multiple time periods so that the resulting decisions may rely on past memories; see \eqref{qinst} for further detail. The NDR stipulates the following form of the control: 
\begin{equation}\label{NDRgeneric}
\mu=\Theta(x,\,\bm q) ,\quad u=\mathcal{P}_{\Omega}\left[\mu \right]
\end{equation}
\noindent where $\Theta$ represents the NDR that maps the states $\bm q$ to the control variables $\mu$; $x$ is the set of parameters of the NDR, which is to be optimized in the off-line training. However, the feasibility of the control $\mu$ in a complex control environment cannot be guaranteed by the NDR, and therefore a projection operator $\mathcal{P}_{\Omega}[\cdot]$ is employed to further map $\mu$ to the feasible control $u$, where $\Omega$ denotes the set of feasible signal control parameters. $\Omega$ may be characterized by fixed cycle time, maximum/minimum green time, and signal offsets, all of which may be expressed linearly. In this case, the projection operator $\mathcal{P}_{\Omega}[\cdot]$ reduces to a quadratic program (see Section \ref{subsecproject} for details).

A NDR of the form \eqref{NDRgeneric} can yield timely signal control decisions given inputs regarding current and past network states, which enables real-time operations as it involves analytical or closed-form transformations. The key step in the NDR approach, which directly impacts its on-line performance, is the optimization of the parameters $x$ through off-line training.

We let $\Phi(\bm q,\,u)=\Phi\big(\bm q, \mathcal{P}_{\Omega}[\Theta(x,\,\bm q] \big)$ be a given network performance measure, which depends on the system state $\bm q$ and the control $u$, along with some inherent uncertainties in the traffic system. For example $\Phi$ may be the delay at a particular junction, or the total emission along a certain corridor. Without loss of generality, we assume that $\Phi$ is subject to minimization.

The problem of optimal NDR can be formulated as
\begin{equation}\label{minphi}
\min_x \Phi\Big(\bm q,\,\mathcal{P}_{\Omega}[\Theta(x,\,\bm q)]\Big)
\end{equation}
\noindent However, note that $\bm q$ is a stochastic variable that varies on a daily basis. For example, $\bm q$ can be the vector of time-varying demands of an arterial network, which vary from day to day. Therefore, a robust feedback control policy such as \eqref{minphi} must take into account the uncertainties in the system. With this in mind, the off-line training of the decision rule may be formulated as the following stochastic optimization problem:
\begin{equation}\label{minphisp}
\min_x \mathbb{E}\big[\Phi\left(\bm q,\,\mathcal{P}_{\Omega}[\Theta(x,\,\bm q)]\right)\big]
\end{equation}
where the objective is to minimize the expectation of the performance measure with uncertain network states $\bm q$.

\subsection{Implementation details}\label{subsecimplementation}

Building on the generic model presented in Section \ref{subsecgenericmodel}, this section presents some implementation details pertaining to the case study of the real-world traffic network in Glasgow presented in  Section \ref{secCaseStudy}.

\subsubsection{Traffic network state variables}
We begin with the state variable $\bm q$, which captures the network-wide traffic state in terms of different measurements (flow, density, speed, etc.) obtained from a network of sensors. Given the discrete time step $t$ ($t$ is an integer) with step size $\delta t$, we express the state variable as
\begin{equation}\label{qinst}
\bm q(t)=
\begin{bmatrix}
q_1(t-n)  &  q_1(t-n+1)  & \ldots   &  q_1(t-m-1)  &  q_1(t-m)
\\
q_2(t-n)  &  q_2(t-n+1)  & \ldots   &  q_2(t-m-1) &  q_2(t-m)
\\
\vdots &  \vdots & \vdots   & \vdots  & \vdots
\\
q_{N}(t-n)  &  q_{N}(t-n+1)  & \ldots   &  q_{N}(t-m-1)  &  q_{N}(t-m)
\end{bmatrix}
\end{equation}
where $0\leq m<n \leq t$ are prescribed integers. On the right hand side of \eqref{qinst} each row corresponds to one sensor, and each column represents one single time step. The integer $n$ is used to indicate the number of past time steps considered when making decisions at the current time step $t$; $m$ is used to account for the fact that data collected in the most recent time intervals may not be immediately available for decision making due to limited capacities of data transmission and computation \citep{Han2017}.  

\begin{remark}\label{rmk1}
In the Glasgow case study presented in Section \ref{secCaseStudy}, the network state is captured by 41 loop detectors, which calculate cross-sectional traffic flows every 2 min (i.e. $N=41$, $\delta t=2$ min). The NDR updates signal timing parameters every 10 min based on the flow information collected in the past 10 min; that is, $m=0$ and $n=4$. 
\end{remark}

\subsubsection{NDR based on feedforward and recurrent neural networks}\label{subsecANN}

\begin{figure}[h!]
	\centering
	\includegraphics[width=\textwidth]{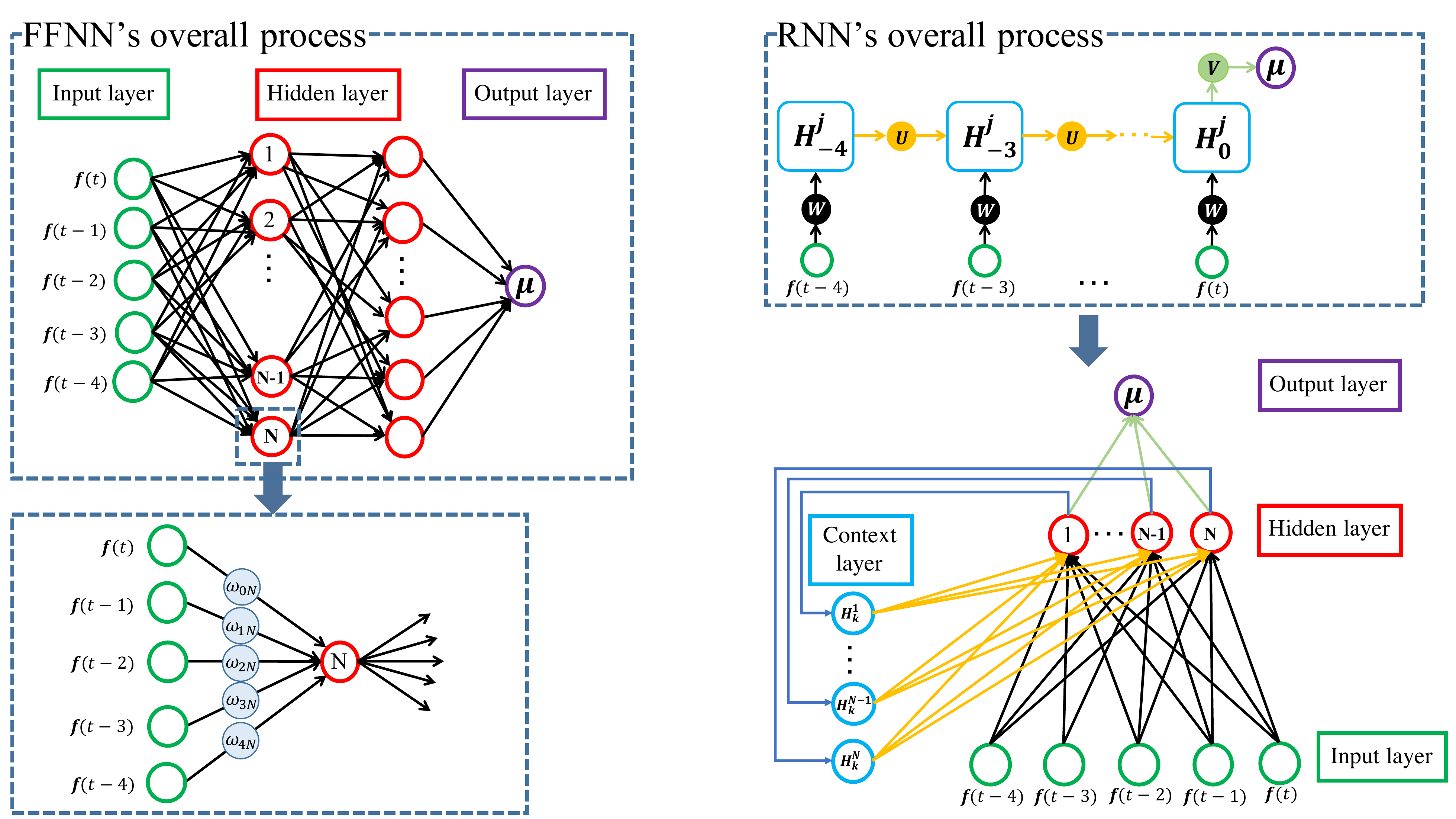}
	\caption{Structures of the FFNN (left) and RNN (right).}
	\label{figNeuralnetwork}
\end{figure}

In this paper, we select feedforward neural network (FFNN) and recurrent neural network (RNN) to instantiate the nonlinear decision rule $\Theta(\cdot\, ,\,\cdot)$. Figure \ref{figNeuralnetwork} illustrates the internal structures of both networks. Given Remark \ref{rmk1}, in order to generate traffic control parameters at time step $t$, both neural networks receive traffic flow vectors in the past 5 consecutive time steps (with a step size of 2 min)
\begin{equation}\label{eqn5f}
\bm f(t),\, \bm f(t-1), \, \bm f(t-2),\, \bm f(t-3),\, \bm f(t-4) \in \mathbb{R}^{41},
\end{equation}
\noindent each being the vector of flows measured at the 41 loop detectors in a 2-min period. \eqref{eqn5f} suggests that the signal control decision made at time $t$ forward depends on the flows in the past five 2-min intervals. To decrease the sensitivity of the neural networks to such input variables, we apply normalization to the vectors $\bm f(t),\, \bm f(t-1)\ldots, \bm f(t-4)$ before feeding them to the neural networks.

The FFNN has two hidden layers with 100 and 50 neurons, respectively. The fully connected neural network employs the Sigmoid activation function, and the weights of connections among the neurons are treated as the parameters $x$ of the NDR in \eqref{NDRgeneric}, to be optimized in the off-line training. In Figure \ref{figNeuralnetwork}(left), the output of the neuron N in the first hidden layer is given as 
\begin{equation}\label{ffnng}
g\left(\sum_{i=0}^4 \omega_{i N}\bm f(t-i) \right)
\end{equation}
\noindent where $\omega_{iN}$'s are the weights, and $g(\cdot)$ is the Sigmoid activation function. Finally, the output $\mu$ for every decision period is generated and used for computing signal control parameters via the projection operator elaborated in Section \ref{subsecproject}.

On the other hand, the RNN has one hidden layer with $N=100$ neurons and one context layer with the same number of neurons as shown in Figure \ref{figNeuralnetwork}(right). The RNN reads the vectors $\bm f(t),\ldots, \bm f(t-4)$ from the input layer one by one in a recursive way:
\begin{equation}\label{eqnrnn}
H_{k}^j=
\begin{cases}
g\left(w_{jk}\bm f(t+k)\right) \quad & k=-4
\\
g\left(w_{jk}\bm f(t+k) +\sum_{i=1}^N u_{ik}H_{k-1}^i \right) \quad & k=-3,\,-2,\,-1,\,0
\end{cases}
\end{equation}
\noindent for $j=1,\,\ldots,\,N$. The RNN iteratively evaluates the quantities $\{H_{k}^j,\, j=1,\ldots, N\}$ after reading a flow vector $\bm f(t-k)$.

In comparison with the FFNN, which perceives all the flow vectors $\bm f(t),\ldots, \bm f(t-4)$ at distinct time steps with a symmetric structure as in \eqref{ffnng}, the RNN processes these vectors in sequence following their chronological order. In this way, the RNN is able to capture the temporal dependencies among these state variables through composition of the activation functions.

\subsubsection{Projection onto the feasible control set}\label{subsecproject}

The signal control parameters typically include cycle time, phasing plans, green times, all-red, and offset \citep{HG2015}. Due to real-world safety considerations reported by \citep{MHHN2015}, the cycle time and phasing plans are fixed in this study. Nevertheless, we note that the NDR framework can be easily extended to dynamically change these control variables.

In this paper, we focus on the phase green times at each and every intersection, denoted $\bm g=(g_1,\, g_2,\,\ldots,\,g_N)^T$, where $N$ is the number of phases. The green times $g_i$ of all the phases must satisfy the following constraints:
\begin{equation}\label{gtcons}
g_{\text{min}}\leq g_i \leq g_{\text{max}} \quad \forall i,\qquad \sum_{i=1}^N g_i= T_{\text{cycle}}- \Delta
\end{equation}
\noindent where $g_{\text{min}}$ and $g_{\text{max}}$ denote minimum and maximum green times, respectively; $T_{\text{cycle}}$ is the fixed cycle time, and $\Delta$ includes amber (all-red) time and pedestrian phase time, which are fixed for safety reasons.

Given the green times $\hat{\bm g}=(\hat g_1,\,\hat g_2,\,\ldots,\,\hat g_N)^T$ as output of the neural network $\Theta(x,\,\bm q)$, which do not necessary satisfy \eqref{gtcons}, the minimum 2-norm projection $\mathcal{P}_{\Omega}$ onto the feasible set can be formulated as the following quadratic program:
\begin{equation}\label{qp}
\min_{\bm g} {1\over 2} \big\| \bm g - \hat{\bm g}\big\|^2= {1\over 2} \big(\bm g-\hat{\bm g} \big)^T \big(\bm g-\hat{\bm g} \big)
\end{equation}
\noindent subject to the linear constraints \eqref{gtcons}. Applying the Karush-Kuhn-Tucker conditions \citep{Friesz2010}, we can explicitly express the solution as:
$$
\bm g=(g_1,\,g_2,\,\ldots,\,g_N)^T , \quad g_i=\big\{\hat g_i -\lambda \big\}_{g_{\text{min}}}^{g_{\text{max}}}
$$
where we employ the notation
$$
\big\{\hat g_i -\lambda \big\}_{g_{\text{min}}}^{g_{\text{max}}}\doteq 
\begin{cases}
g_{\text{min}}  \quad & \hbox{if} \quad  \hat g_i -\lambda < g_{\text{min}}
\\
\hat g_i -\lambda \quad & \hbox{if} \quad  g_{\text{min}}\leq \hat g_i -\lambda \leq g_{\text{max}}
\\
g_{\text{max}}  \quad & \hbox{if} \quad  \hat g_i -\lambda > g_{\text{max}}
\end{cases}
$$
and the dual variable $\lambda$ is such that
\begin{equation}\label{algeqn}
\sum_{i=1}^N \big\{ \hat g_i -\lambda \big\}_{g_{\text{min}}}^{g_{\text{max}}}= T_{\text{cycle}}-\Delta,
\end{equation}
which can be found by numerically solving the algebraic equation \eqref{algeqn}. Note that in reality the maximum and minimum green times may vary across different signal phases, in which case the formulae above remain valid.

\subsection{Off-line optimization of the NDR}\label{subsecofflinetraining}
This section presents details of the simulation-based optimization procedure, which serves as the off-line module to train and optimize the NDR; i.e. the neural network presented in Section \ref{subsecANN}. The main purpose is to find the optimal (or near optimal) solutions of the optimization problem \eqref{minphisp}, which is recapped here:
\begin{equation}\label{minphispcopy}
\min_x \mathbb{E}\big[\Phi\left(\bm q,\,\mathcal{P}_{\Omega}[\Theta(x,\,\bm q)]\right)\big]
\end{equation}
The inherent stochasticity in the network states $\bm q$ can be handled in different ways such as using robust optimization and stochastic optimization, with varying degrees of conservatism and computational complexity; see \cite{BBC2011, LHGFY2015} and \cite{Han2017} for more discussions. In this paper, due to the potentially expensive evaluation procedure, which is done through microscopic traffic and emission simulations, we propose a Monte-Carlo type evaluation method.

Specifically, the overall optimization procedure, which is viewed as the off-line module of the proposed signal control framework, can be divided into two levels; see Figure \ref{figworkflow}. The upper-level problem is to find the optimal parameters $x$ to minimize the expectation shown in \eqref{minphispcopy}. The objective function involves traffic micro-simulation and high-fidelity emission modeling, whose dynamics and uncertainties are difficult to characterize analytically. Therefore, we employ a heuristic method based on Particle Swarm Optimization method (PSO) to find optimal $x$. The PSO is chosen here as it requires only zeroth-order information of the objective and the constraints. In addition, although the performance of PSO varies depending on the application or parameters, research shows evidences of PSO or its variants outperforming other metaheuristic or evolutionary algorithms such as ant colony optimization, simulated annealing, tabu search, and genetic algorithm \citep{Yin2006, SRV2010, SH2008}. On the other hand, the lower-level problem seeks to evaluate the expected network performance (in terms of traffic and emission indicators) with given parameters $x$, while taking into account stochasticity in the traffic states $\bm q$ and microscopic traffic dynamics such as driving behavior and route choices.

\subsubsection{Particle Swarm Optimization}\label{subsecPSO}

Particle Swarm Optimization (PSO) \citep{BVA2007} offers an efficient and flexible trade-off between optimality of the solution and computational resources. It is based on the social behavior of a group of animals, called a {\it swarm}. In a swarm, the animals are represented as particles, and can collaborate and share information to adjust their positions in the search for a certain location. The adjustment of their positions is based on the swarm's collective memory on the best location attained so far (hereafter referred to as ``gbest"), and the individual memory of the best location that the particle has attained so far (hereafter referred to as ``pbest"). As a result of the position adjustment, the particles tend to converge to either $G$ or $P_j$. Although the performance of PSO varies depending on the domain of applications or parameters chosen, research shows evidence of PSO or its variations outperforming well-established metaheuristics such as genetic algorithm, ant colony optimization, simulated annealing, and tabu search.

Given the objective function to be minimized, denoted $f(\cdot)$, and the feasible domain $S$, the following pseudo code summarizes the PSO procedure.

\begin{framed}
  \noindent {\bf Particle Swarm Optimization}
  \begin{description}
  \item[\bf Input.] Population size $N$, $\{\omega_k:\, k\geq 0\}\subset (0,\,1)$, $c_1,\,c_2>0$.

  \item[\bf Step 0.] Let $k=0$. Randomly initialize the particles'
    positions $X_i^0$ and velocities $V_i^0$, $1\leq i\leq N$. Initialize ``pbest" $P_i^0$ and ``gbest" $G^0$ as follows:
    \begin{equation*}
      P_i^0=X_i^0 \quad 1\leq i\leq N,  \qquad G^0=P_{i^*}^0
    \end{equation*}
    \noindent where $i^*=\underset{1\le i\le N}{\text{argmin}}\,f(P_i^0)$.

  \item[\bf Step 1.] Update the velocities and positions: for all $1\le i\le N$, 
    \begin{align*}
      V_{i}^{k+1}=&\omega_kV_{i}^{k}+c_1 r_1(P_i^k-X_i^k)+c_2r_2(G^k-X_i^k)
      \\
      X_i^{k+1}=&\mathcal{P}_{S}[X_i^k+V_i^{k+1}] 
    \end{align*}
    where $r_1$ and $r_2$ are random numbers uniformly generated within
    $[0,\,1]$.

  \item[\bf Step 2.] Evaluate the objective values $f(X^{k+1}_i)$
    for all $1\leq i\leq N$.

  \item[\bf Step 3.]  Update ``pbest" and ``gbest":
    \begin{align*}
      P_i^{k+1}=&
      \begin{cases} 
        X_i^{k+1}   \quad & \hbox{if} ~ f(X^{k+1}_i)< f(P_i^{k})
        \\ 
        P_i^{k}  \quad & \hbox{Otherwise}
      \end{cases}\qquad\quad \forall 1\leq i\leq N  \nonumber
      \\
      G^{k+1}=&
      \begin{cases} 
        P_{i^*}^{k+1} \quad  &  \hbox{if}~ \underset{1\leq i\leq N}{\min}
        f(P_{i}^{k+1}) <  f(G^{k})
        \\
        G^{k}    \quad & \hbox{Otherwise}
      \end{cases}
    \end{align*}
    where $i^*=\underset{1\leq i\leq N}{\text{argmin}}\, f(P_{i}^{k+1})$

  \item[\bf Step 4.] If the stopping criterion is met (e.g. no improvement
    in the objective within a given number of consecutive iterations),
    terminate the algorithm with output $G^{k+1}$. Otherwise, let $k=k+1$,
    and go to Step 1.
  \end{description}
\end{framed}

\subsection{Off-line training procedure}

\begin{figure}[h!]
\centering
\includegraphics[width=0.8\textwidth]{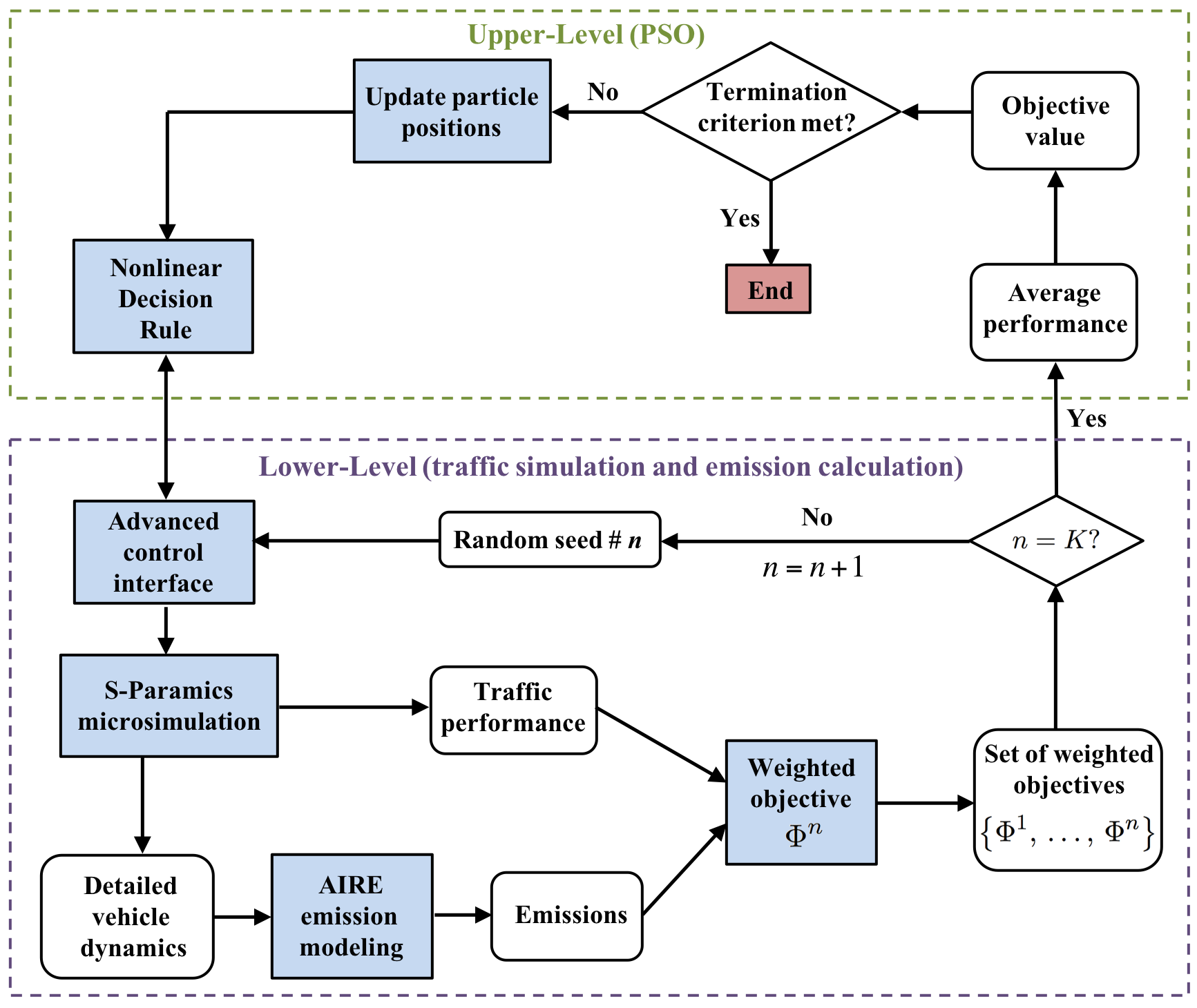}
\caption{Off-line training (optimization) procedure of the nonlinear decision rule.}
\label{figworkflow}
\end{figure}

The off-line training of the NDR amounts to a simulation-based optimization procedure, which requires PSO to be carried out in conjunction with the Monte-Carlo approach that assesses the NDR with given parameters (for FFNN or RNN) via microsimulation and emission calculation. The work flow of the simulation-based optimization is outlined in Figure \ref{figworkflow}, with individual key components explained below.

\subsubsection{PSO for solving optimization \eqref{minphisp}} 

The PSO is an agent-based search method, which is detailed in Section \ref{subsecPSO}. In each iteration of the PSO, a total of $N$ agents conduct independent search by evaluating, for a given NDR, the corresponding objective value, which is defined as the expectation in Equation \eqref{minphisp}. The stochasticity arises from the microscopic traffic simulation where the departure rates and route choices are randomly sampled based on an origin-destination matrix describing travel demands. Another source of stochasticity comes from the microscopic driving dynamics, which involve car-following, lane-changing, and gap-acceptance behavior. All these stochasticities in each simulation run are populated by a random seed, and we use $K$ distinct random seeds to represent the stochastic nature of the traffic states. The aforementioned expectation is then approximated as the average over $K$ independent simulation runs.

 In the case study presented in Section \ref{secCaseStudy}, the PSO employs a population size of $N=5$, and the algorithm is terminated if no improvement is made on the objective within 20 iterations or when the total iteration number reaches 45.

\subsubsection{Traffic simulation}

 The microscopic traffic simulation is performed using the S-Paramics software \citep{Paramics2011}, which not only calculates various traffic {\it key performance indicators} (KPIs) such as travel time, delay and throughput, but also produces detailed vehicle trajectories at a resolution of 0.5 second, which are used as input of emission modeling. For the case study presented in the next section, a microsimulation model is set up for the west end of Glasgow, and calibrated using a combination of macroscopic and microscopic data. See Section \ref{subsecsimulationglasgow} for more details.

The number of traffic simulation (and emission estimation) that needs to be performed within one major PSO iteration is equal to $N\times K$ where $N$ is the population size (independent search agents) and $K$ is the number of random seeds used to populate stochastic parameters and dynamics in the simulation.

\subsubsection{Emission calculation}
A main feature of the proposed real-time signal control framework is the consideration of environmental impact caused by exhaust emissions from vehicles, which is directly impacted by vehicle dynamics and the signal control strategies. In this paper, we focus on CO$_2$ and Black Carbon (BC). CO$_2$ is the primary greenhouse gas and contributes to global warming, while BC causes serious health concerns such as respiratory problems, heart attacks and lung caners. 

We use the AIRE (Analysis of Instantaneous Road Emissions) vehicle exhaust emission model \citep{AIRE2011} to calculate the instantaneous total carbon and particulate matter emissions resulting from the combustion of fuel throughout each vehicle journey in the simulation. AIRE is a subsidiary software program that interfaces with S-Paramics and post-processes the output of the traffic simulation. Through built-in Instantaneous Emissions Modeling (IEM) tables, AIRE is able to generate estimated value of vehicle emission for each simulated vehicle \citep{TS2011}.

As AIRE does not calculates CO$_2$ and BC directly, we further develop a post-processing tool to convert total carbon and particulate matters into CO$_2$ and BC emissions. In particular, the following procedures are followed to achieve this.
\begin{itemize}
\item The total carbon metric is based on the PHEM (Passenger Car and Heavy Duty Emission Model) fuel consumption metric and consequently can be directly converted into a representative CO$_2$ emissions.  This is done by using the atomic weights of Carbon and Oxygen to generate a factor of 44/12 (one molecule CO$_2$ weighs 44, one atom carbon weighs 12). 
\item The calculation of BC is based on the estimated PM$_{10}$ emission rates using the COPERT IV methodology for conversion \citep{GKNS2006, MHHNVTBL2016}.
\end{itemize}

\subsubsection{Weighted objective}

To simultaneously reduce traffic congestion and emissions, we reformulate the multi-objective optimization problem into a single-objective one through  scalarization:
\begin{equation}\label{objscalar}
\min \text{Delay},\quad \hbox{or}\quad \min w_1\cdot {\text{Delay}\over n_1} + w_2 \cdot{\text{CO}_2 \over n_2} \quad \hbox{or}\quad  \min w_1\cdot {\text{Delay}\over n_1} + w_3\cdot {\text{BC}\over n_3} 
\end{equation}
where $\text{Delay}$ refers to network-wide average delay per vehicle, CO$_2$ and BC are respectively the network-wide total CO$_2$ and BC emissions. Constants $n_1$ and $n_2$ are normalization factors to bring the three objective values to a comparable numerical scale. $w_i$'s are positive weights.

\subsubsection{Advanced control interface (ACI)}
ACI is a method of accessing the traffic model via external program and exchanging information. In S-Paramics, it uses a component protocol called SNMP (Simple Network Management Protocol) to achieve that. Through this protocol, external program can organize, collect, and modify traffic information to change the condition of traffic model \citep{Paramics2011}. For example, in this paper, the ACI has two main functions; parameter/data exchange and simulation synchronization. First, the program will access real-time (in simulation) traffic data to monitor the performance of the traffic network. This information will be saved and used for the responsive signal optimization procedure. Second, the information will be sent to S-paramics for the synchronization purpose, as the right signal timings need to be implemented in the right time during the simulation. The ACI model has been developed by using visual basic application (VBA), which is an integral part of our experiment set up and facilitates the information exchange and control among different models.

\section{Case Study in Glasgow}\label{secCaseStudy}

\subsection{Simulation of the test site}\label{subsecsimulationglasgow}

\begin{figure}[h!]
\centering
\includegraphics[width=.88\textwidth]{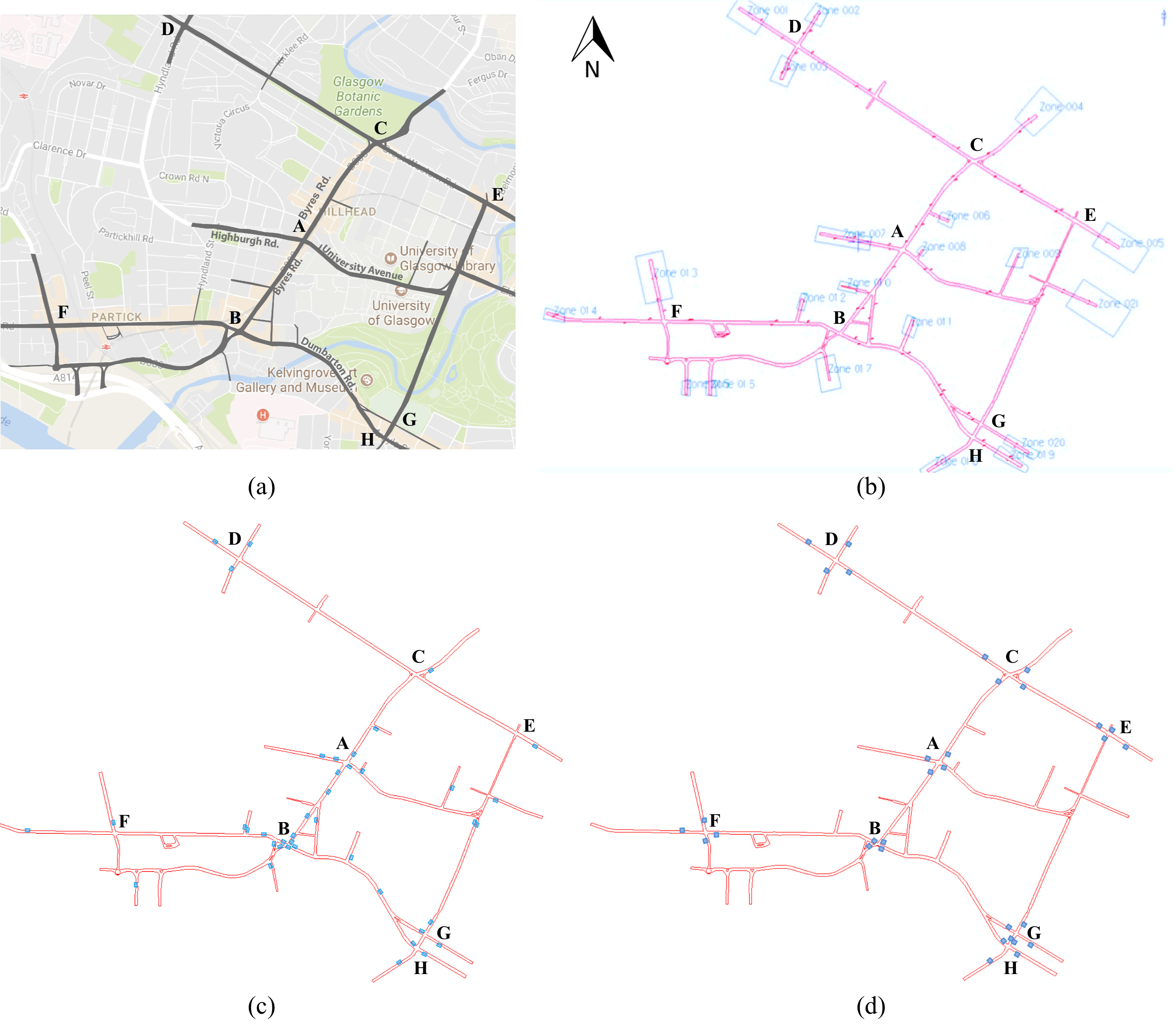}
\caption{(a) The test area in Glasgow with 8 signalized intersections. (b) Road network with 21 Zones. (c) The locations of the 41 loop detectors in the real world. (d) Alternative locations of the loop detectors for the comparative study.}
\label{figGlasgownetwork}
\end{figure}

The proposed real-time signal control framework has been applied to a real-world test network in Glasgow, Scotland. The traffic simulation modeling was conducted within the EU-funded CARBOTRAF project, which aims to support adaptive traffic management for reducing urban congestion and associated environmental and health impacts. The study area is the west part of Glasgow (see Figure \ref{figGlasgownetwork}) with 14 signalized junctions and 478 links. There are 21 zones (Figure \ref{figGlasgownetwork}(a)), giving rise to 420 origin-destination pairs.

\begin{figure}[h!]
\centering
\includegraphics[width=\textwidth]{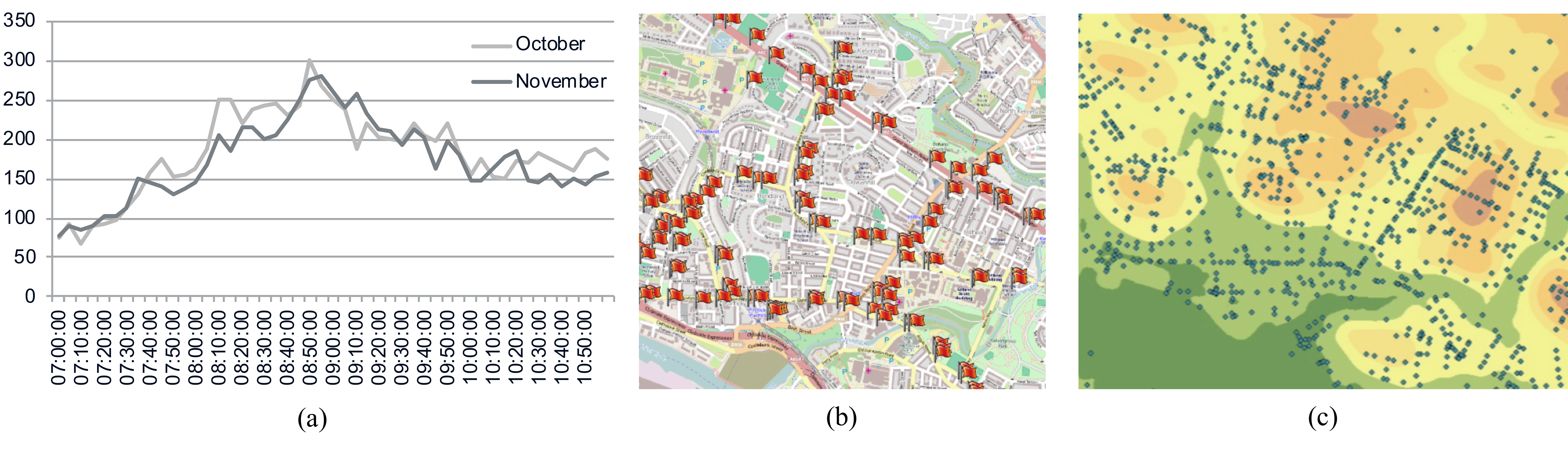}
\caption{(a): 5-min average traffic flow (veh/hr) on Byres Road. (b): Bus stops in and around the test network. (b): Network nodes overlaid with Digital Elevation Model.}
\label{figsite}
\end{figure}

A typical demand scenario for the test network was generated within the S-Paramics simulation software for 7:30-9:30am, which represents morning peak of a typical working day (Monday-Thursday) in 2010; see Figure \ref{figsite}(a). The microscopic model has been built using the OS-ITN network to represent the supply, and a seeded demand matrix obtained from loop detectors that represent the within-day and day-to-day variability of traffic \citep{MHHN2015}. For the baseline control scenario, we consider the default traffic signal timing plans provided by the Glasgow City Council. The baseline model has been fully calibrated and validated using a combination of loop detector and floating car data \citep{MHHN2015}.

The vehicle fleet that has been simulated consists of private cars, taxis, buses, vans, light goods vehicles, and heavy goods vehicles. The fleet composition is based on the Annual Average Daily Flow (AADF) data for Glasgow city from the Department for Transport between 2000 and 2010. This allows us to capture realistic traffic dynamics with mixed vehicle types and to accurately estimate vehicle emissions with detailed emission factors for different vehicle types. Moreover, road gradient has been explicitly modeled based on the Digital Elevation Model as it has been shown to play a significant role in engine load and, subsequently, carbon emissions \citep{SMH2016}. Figure \ref{figsite} shows the bus stops in and around the test network as well as the digital elevation information for the study area.

\subsection{Signal control details}

\begin{figure}[h!]
\centering
\includegraphics[width=\textwidth]{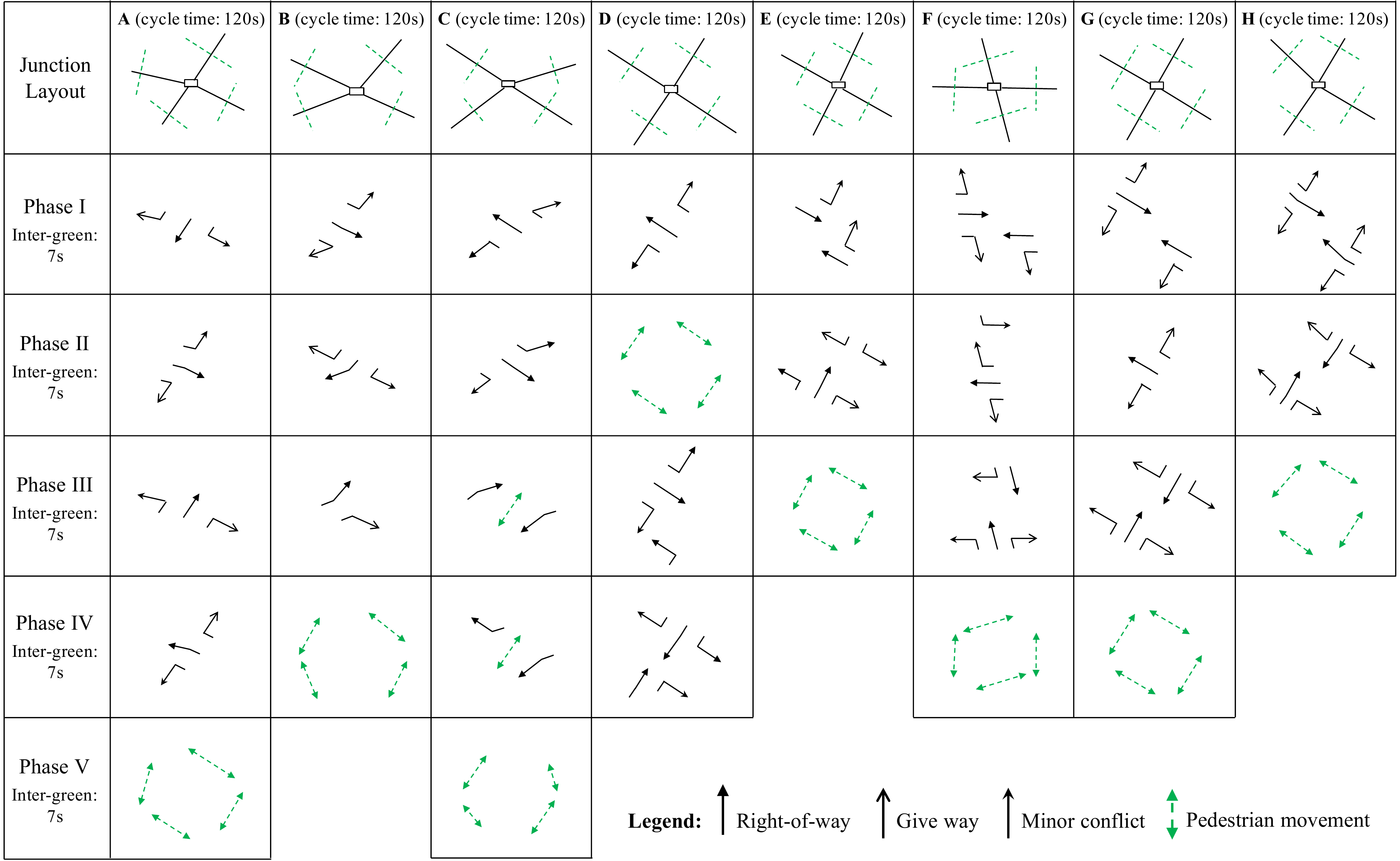}
\caption{Phasing plans of the eight signal intersections.}
\label{figphasing}
\end{figure}

The network has eight major signal intersections, shown as intersections A-H in Figure \ref{figGlasgownetwork}. Other minor junctions in the network are either priority junctions/roundabouts or controlled by actuated signals, which are excluded from our control framework. The cycle times, inter-green (including amber and all-red), and phasing schemes of the eight main intersections are shown in Figure \ref{figphasing}. These quantities are fixed in our NDR framework per real-world control and safety requirements imposed by the Glasgow City Council (GCC), and parameters subject to real-time optimization are the green times of all the vehicle-movement phases. Note that signal offsets could be easily included as additional decision variables in our control framework, but they are difficult to be adjusted dynamically within the microsimulation. For this reason the offsets are all fixed using the default setting (provided by the GCC) in our control framework.

The resolution for the adaptive signal control is 10 min, which means that the signal timings are adjusted every 10 minutes depending on the real-time traffic conditions. Accordingly, 10-min average traffic flow information collected by the 41 loop detectors (see Figure \ref{figGlasgownetwork}(b)) is provided to the NDR to update signal controls for the next 10 min.

\subsection{Signal control scenarios}\label{subsectestscenarios}

To investigate the extent of traffic and environmental impact of the proposed real-time signal controls, and to make a case for coordinated signal controls on a network-wide level, we consider three test scenarios with varying controllability.

\begin{itemize}
\item[(1)] {\bf Junction Level [JL]}: only intersection A is controlled dynamically by our NDR approach; all the other seven intersections are controlled by the default signal timings (provided by the GCC). Intersection A is of critical importance as it connects traffic from the west to major local destinations including universities and hospitals. In the real-world, location A is most affected by traffic congestion and air pollution.

\item[(2)] {\bf Corridor Level [CL]}: only intersections A, B and C are dynamically controlled by our NDR approach; the other five intersections are controlled by the default timings. Intersections A-C are located along the Byres road, which is a strategic corridor connecting the radial routes to the city center for drivers approaching from the west. 
 
\item[(3)] {\bf Network Level [NL]}: all eight junctions in the network are simultaneously and dynamically controlled by the NDR approach. In this way, the signal timings are coordinated in a centralized fashion. In other words, the control of any local intersection is informed by the traffic states and other signal timing plans on the entire network. This is in contrast to distributed controls, which seeks local efficiency over global optimality. 
\end{itemize}

\noindent As a benchmark for comparing with the proposed signal control strategy, we consider:
\begin{itemize}
\item[(4)] {\bf Glasgow City Council [GCC]}: the fixed signal timing plan provided by the Glasgow City Council, which is derived from static OD route flow information. 
\end{itemize}

In accordance with the aforementioned control scenarios, we conduct off-line training (optimization) of the NDR by minimizing the average delay, CO$_2$ emission, BC emission through a weighted combination of these objectives. This allows us to understand the potential trade-off between traffic efficiency and environmental impact. Then, the on-line performance of the optimized NDR is tested in 30 independent simulation runs with 30 random seeds that are different from the ones used in the training. The resulting performance of the traffic network, measured in terms of delay, CO$_2$ emissions, BC emissions, queuing and throughput, is presented in the following sections.

\subsection{Test results and discussion}
The test results are evaluated against four key performance indicators (KPIs):
\begin{itemize}
\item network-wide average delay. The delay is defined as the difference between the actual journey time of a trip minus the free-flow time obtained by assuming little traffic; 
\item network throughput, defined to be the number of vehicles completing their trips by the end of the simulation period;
\item average vehicles in queue, which is defined on a link level; and
\item network-wide CO$_2$ and black carbon emissions. 
\end{itemize} 

\subsubsection{Overall performance of the proposed signal controls}

\begin{figure}[h!]
	\centering
	\includegraphics[width=.9\textwidth]{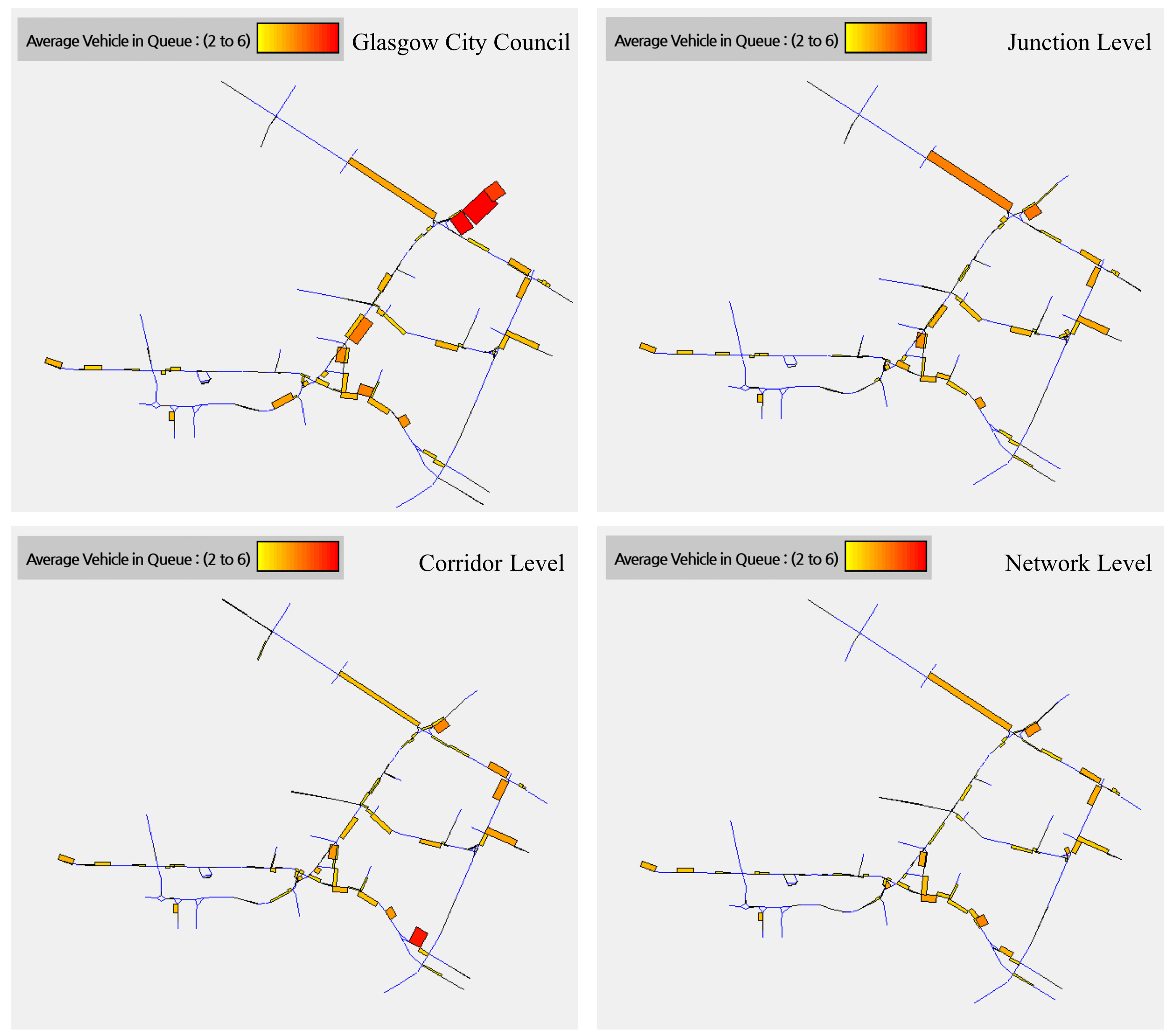}
	\caption{Average number of vehicles in queue.}
	\label{figaveragevehicleinqueue}
\end{figure}

Figure \ref{figaveragevehicleinqueue} shows the average number of vehicles in queue on each link of the network, which is a direct indicator of network congestion. In the case of [GCC], significant congestion is seen along the Byres corridor, especially on the northern entrance. For the proposed methods, widening the scope of the signal controls ([JL] to [CL] to [NL]) tend to mitigate the congestion on the network level overall. However, minor spatial trade-offs of congestion can be seen, for example, between [CL] and [JL]. Through a coordinated control of the three intersections A, B and C, [CL] effectively reduces the congestion on the Byres corridor, especially on the northern entrance (including the Great Western Rd.) compared to [JL]. However, more significant queueing on the southeast part of the network results from the [CL], possibly due to (1) lack of direct control of that area; and (2) increased traffic flow on the Dumbarton Rd. as a result of improved Byres corridor. Such trade-off of congestion at different parts of the network reveals the complexity of network-wide adaptive signal control as drivers' route choices are affected by real-time traffic conditions \citep{HSLFY2015}. Finally, [NL] eliminates all the major queuing on the network and achieves the highest efficiency in terms of vehicle queues. Nevertheless, even in this case some queuing still remains along the northern corridor (Great Western Rd.); this is due to the lack of sufficient sensor coverage along this main corridor; see Figure \ref{figGlasgownetwork}(b), and the proposed signal controls are not fully informed by the traffic states there.

\begin{figure}[h!]
	\centering
	\includegraphics[width=\textwidth]{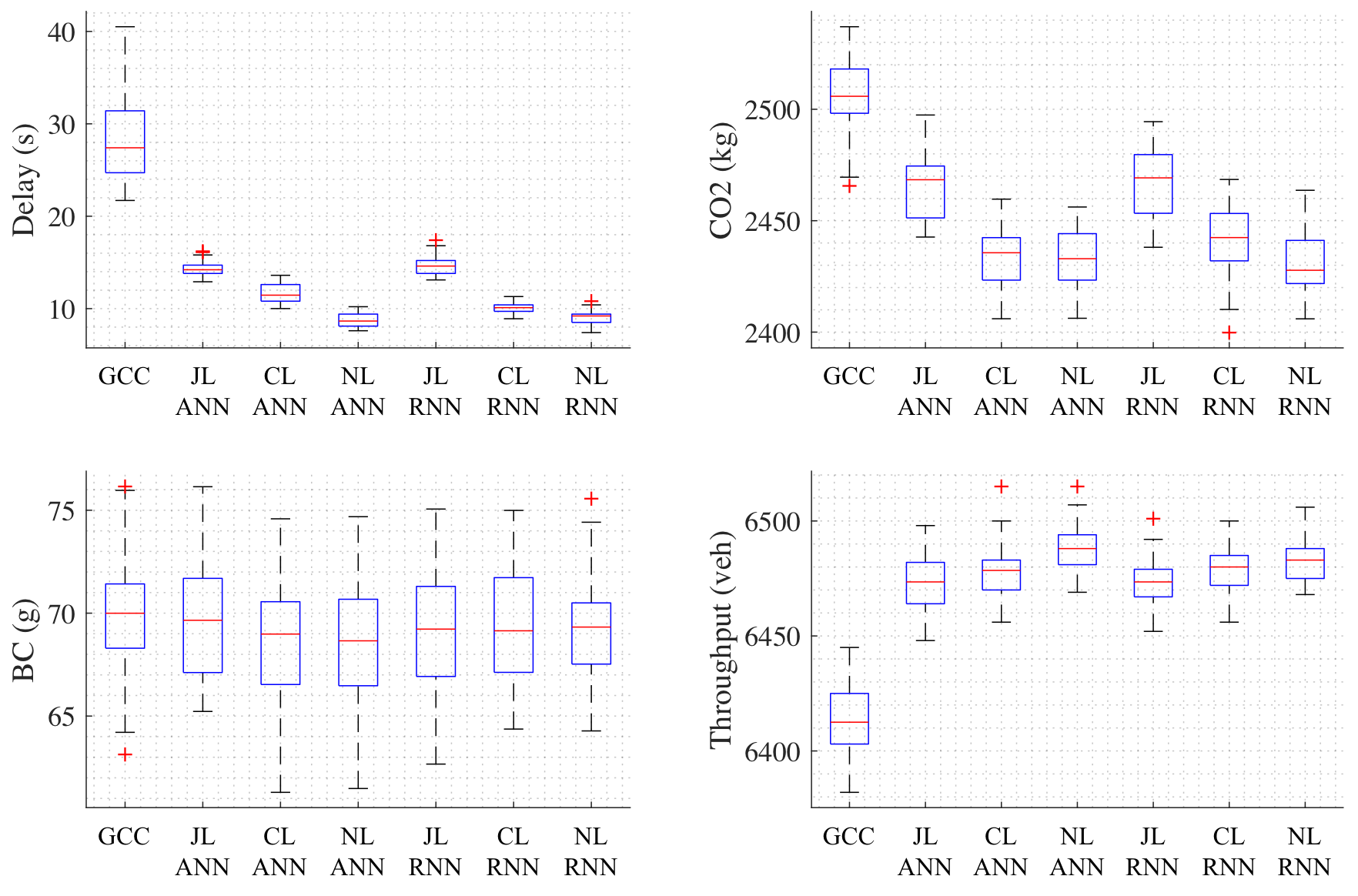}
	\caption{Box plot summary (with 30 random simulation runs) of the performance of the four control scenarios in terms of average network delay, total carbon emission, black carbon emission, and throughput.}
	\label{figBoxplots}
\end{figure}

\newpage
Figure \ref{figBoxplots} shows the performances of the four control scenarios (GCC, JL, CL, NL) based on FFNN and RNN, in terms of delay, throughput, total carbon and black carbon emissions, followed by Table \ref{tab:improvementtable} summarizing the average improvements of the proposed signal controls over the baseline scenario (GCC). It can be seen that the proposed signal control methods significantly outperform the existing signal control (GCC). Among all four KPIs vehicle delay has the most drastic improvement, from around 28 seconds per vehicle to below 10 seconds (NL). This is followed by CO$_2$ emission and network throughput, with up to 73 kg reduction and 74 veh increase, respectively. The decrease in CO$_2$ is likely caused by increased travel speeds as a result of reduced congestion, as CO$_2$ emissions tend to increase at low driving speeds \citep{LEFE2011}. The decrease of black carbon is comparatively less significant with 0.5-1.4 g reduction. Black carbon forms during incomplete combustion of carbonaceous fuels, and is primarily caused by sudden acceleration and brake of vehicle movements  \citep{MHHNVTBL2016}. Therefore, the reduction of BC is more significant at local intersections than on the network level (see Figure \ref{figjunctionemission}). It is also clear from Figure \ref{figBoxplots} and Table \ref{tab:improvementtable} that the benefits of the proposed real-time signal control are pronounced when more signalized intersections (e.g. network-wide control with 8 signals) are simultaneously controlled.

\begin{table}[h!]
	\centering
		\caption{Statistical summary (with 30 random simulation runs) of the improvement of [JL], [CL] and [NL] over baseline [GCC]}
		\begin{tabular}{| c   c | l | l | l |}
			\hline
			 &   & Scenario [JL] & Scenario [CL] & Scenario [NL]  
			 \\\hline
\multirow{2}{*}{Delay}  & \multicolumn{1}{|c|}{FFNN}  & 14.0 s (48.2\%) & 16.7 s (58.1\%) & 19.6 s (68.4\%)  
\\
\multirow{1}{*}{}           &  \multicolumn{1}{|c|}{RNN} & 13.7 s (47.2\%) & 18.3 s (63.7\%) & 19.3 s (67.2\%)  

\\\hline

\multirow{2}{*}{CO$_2$}   & \multicolumn{1}{|c|}{FFNN} & 38.5 kg (1.5\%) & 69.9 kg (2.8\%) & 71.3 kg (2.8\%)  
\\
\multirow{1}{*}{}                & \multicolumn{1}{|c|}{RNN} & 35.3 kg (1.4\%) & 63.8 kg (2.5\%) & 73.2 kg (2.9\%)  

\\\hline

\multirow{2}{*}{BC}           &\multicolumn{1}{|c|}{FFNN}    & 0.59 g (0.7\%) & 1.3 g (1.8\%) & 1.4 g (2.0\%)  
\\
\multirow{1}{*}{}                &\multicolumn{1}{|c|}{RNN}    & 0.81 g (1.2\%) & 0.72 g (0.8\%) & 0.72 g (0.8\%)  

\\\hline

\multirow{2}{*}{Throughput}  &  \multicolumn{1}{|c|}{FFNN} &  59.6 veh (0.9\%) & 64.7 veh (1.0\%) & 73.8 veh (1.2\%)  
\\
\multirow{1}{*}{}                    &  \multicolumn{1}{|c|}{RNN} &  59.8 veh (0.9\%) & 65.1 veh (1.0\%) & 68.2 veh (1.1\%)  

\\\hline
		\end{tabular}
		\label{tab:improvementtable}
\end{table}

\subsubsection{Improvement at junction level}

\begin{figure}[h!]
	\centering
	\includegraphics[width=.6\textwidth]{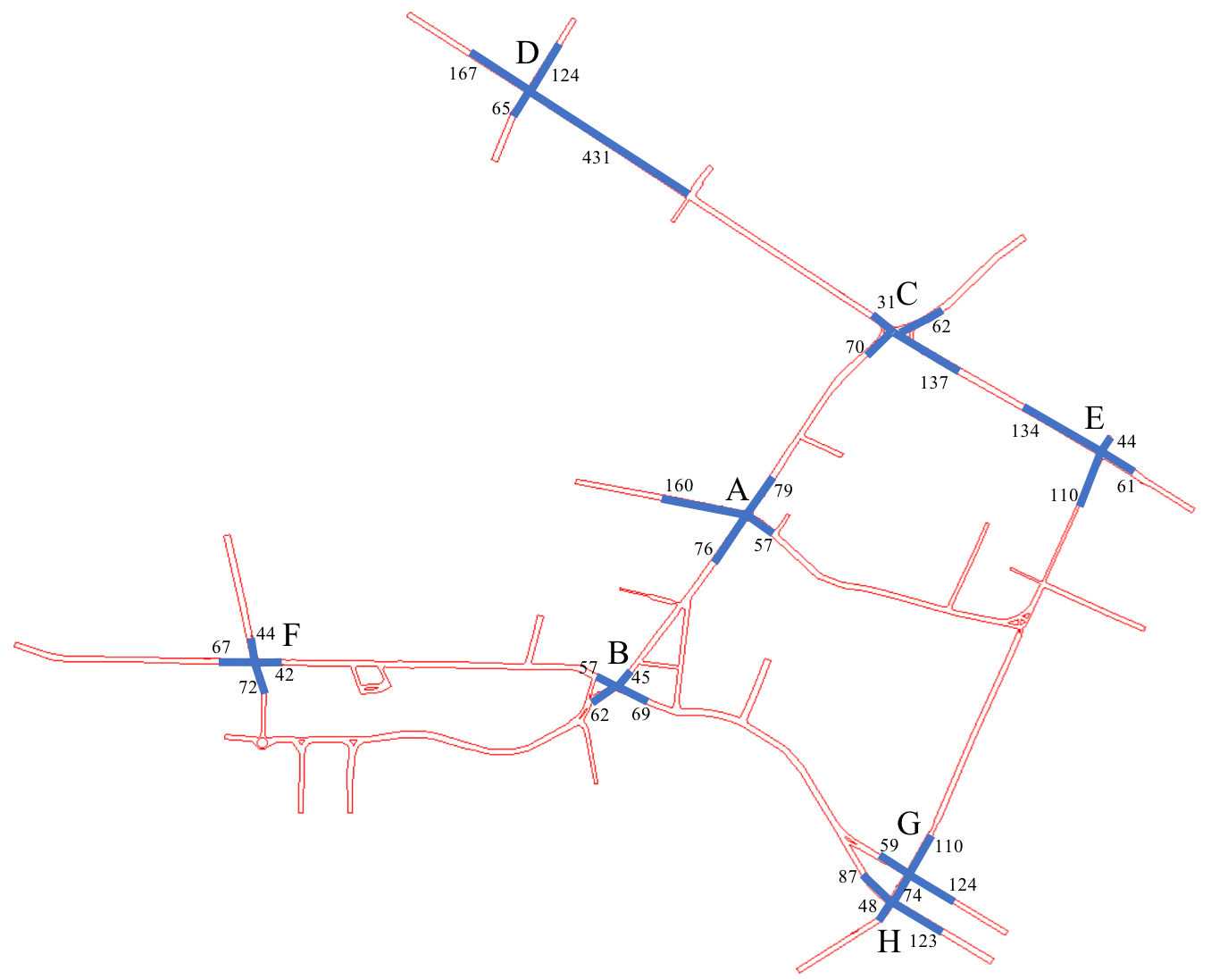}
	\caption{The emissions at signalized intersections are calculated on the highlighted incoming approaches, with their lengths shown (in meter).}
	\label{figsignalrange}
\end{figure}

To further examine the effects of the proposed controls, we evaluate the emission reduction at individual signalized intersections. The emission at an intersection is calculated as the sum of emissions at its incoming approaches, as shown in Figure \ref{figsignalrange}.

\begin{figure}[p!]
	\centering
	\includegraphics[width=\textwidth]{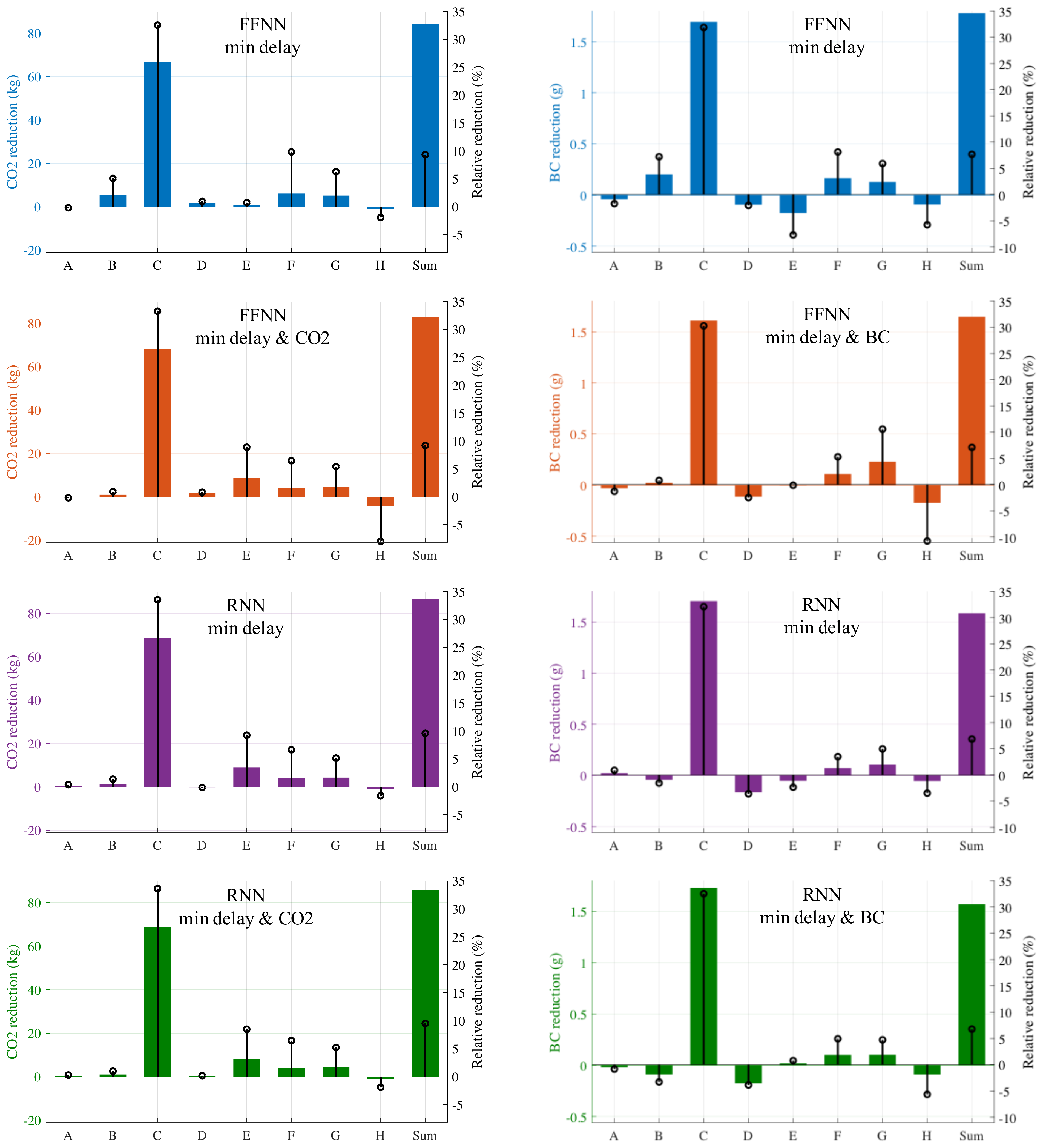}
	\caption{Reductions of CO$_2$ and BC emissions at individual intersections.}
	\label{figjunctionemission}
\end{figure}

In Figure \ref{figjunctionemission}, we show the average absolute (left axis) and relative (right axis, in \%) CO$_2$ and BC reductions at the eight signal intersections. The cases examined include: FFNN vs. RNN, and different objectives in the off-line training (i.e. optimizing delay only or combination of delay and CO$_2$/BC). Both absolute and relative CO$_2$/BC reductions at the junction level show far greater improvement compared to the network level (see Table \ref{tab:improvementtable} for comparison). The overall reduction of CO$2$ (BC) at the junctions is above 80 kg (1.5 g), when the network-level reductions are up to 73 kg (1.4 g). This means that the network-wide reduction of emissions is almost entirely attributed to the improved signal controls at individual intersections, which indeed shows the effectiveness of the proposed controls. In addition, the majority of the savings occur at junction C, with over 30\% reduction of both CO$_2$ and BC. As for rest of the intersections, the reductions of CO$_2$ are mostly positive except H, while the reductions of BC are mixed. Finally, in terms of the optimization objective, minimizing delay alone seems to yield similar emission reductions as the weighted sum of delay and emissions.

\subsubsection{Trade-off between vehicle delay and emissions}

\begin{figure}[h!]
	\centering
	\includegraphics[width=.8\textwidth]{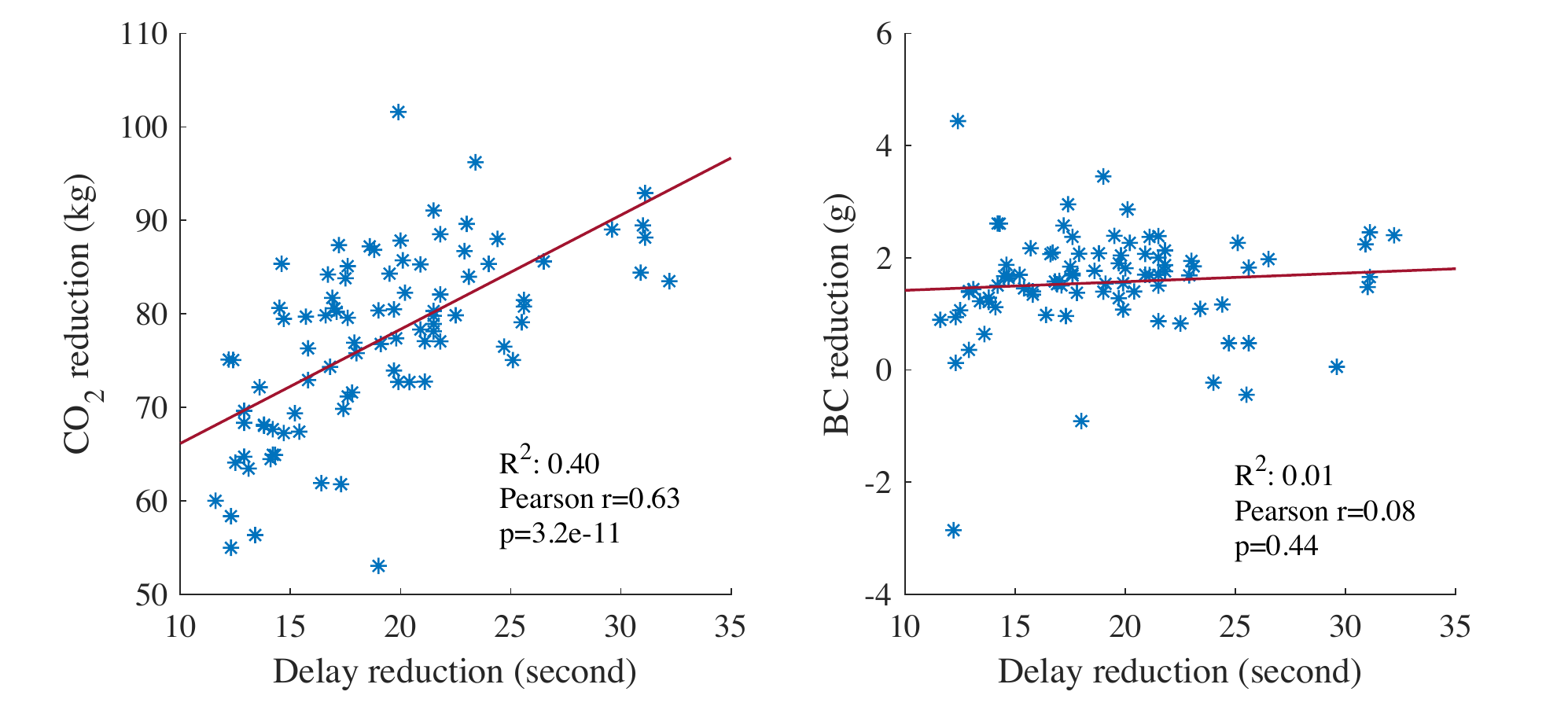}
	\caption{Correlation analysis of delay reductions vs. emission reductions at the junction level.}
	\label{figtradeoff}
\end{figure}

As Figure \ref{figjunctionemission} suggests, minimizing delay as the only objective seems to yield similar levels of emission reduction than the joint minimization of delay and emission. Intuitively, reducing vehicle delays leads to increased average speed, which could reduce CO$_2$ emissions, and reduces vehicle idling and acceleration/deceleration events. To further investigate the potential correlation between vehicle delays and emissions, in Figure \ref{figtradeoff} we show the scatter plots of delay reduction vs. total emission reduction at the junction level. These data points are obtained from a total of 90 independent on-line simulation runs, where the NDR was respectively optimized off line with the three objectives shown in \eqref{objscalar}. The figure shows that reductions of CO$_2$ are positively correlated with delay reductions, as indicated by the Pearson test ($p\approx 0$). This is consistent with the interpretation that CO$_2$ emissions are dependent on average vehicle speed, which is related to vehicle delays. On the other hand, BC reductions do not show meaningful correlation with delay reductions ($p=0.44$). This is attributed to the fact that BC emissions are primarily caused by stop-and-go cycles and highly dependent on vehicle fleet composition (e.g. buses, HGVs), which are not directly related to average vehicle delays.

We also observe from Figure \ref{figjunctionemission} that, when jointly minimizing delay and emissions (CO$_2$ or BC), the reduction of emissions does not improve compared to the case when only delay is minimized. Aside from the possibility that the PSO-based off-line heuristic optimization does not yield global optimal within the required computational resources, the lack of discernible trade-off (i.e. statistically significant negative correlation) between delay reduction and emission reduction, as seen in Figure \ref{figtradeoff}, also suggests that optimizing delays in this case seems to be sufficient in reducing emissions at signalized intersections. More effective measures for emission reduction may involve localized, actuated controls such as transit signal priority and offset optimization, which are beyond the scope of this paper.

\subsubsection{Effect of sensor locations}

The proposed NDR framework for real-time signal control perceives the traffic state via the input vector $\bm q$ in \eqref{NDRgeneric}, which, in our case, represents the traffic flows at several key locations in the past 10 minutes. All the results presented to this point are based on the real-world configuration of 41 loop detectors as shown in Figure \ref{figGlasgownetwork}(c). To assess the impact of sensor locations on the performance of the NDR approach, we consider an alternative configuration in Figure \ref{figGlasgownetwork}(d), where each incoming approach of every signalized intersection has a loop detector, and there are 32 such detectors. The real-world sensors are distributed unevenly across the network, with missing information on incoming traffic at several key intersections (C, D, E, F). The alternative sensor configuration makes sure that the controller receives traffic information at all relevant incoming approaches.

\begin{table}[h!]
	\centering
		\caption{On-line performances of the NDR based on real-world and alternative sensor locations (respectively (c) and (d) in Figure \ref{figGlasgownetwork}).}
\begin{tabular}{| c   c | l | l | l | l |}
			\hline
			 &   & Delay & CO$_2$ & BC   & Throughput
			 \\\hline
\multirow{2}{*}{Real-world}  & \multicolumn{1}{|c|}{FFNN}  & 8.7 s & 2,432 kg& 68.6  g& 6,487 veh
\\
\multirow{1}{*}{}                    &  \multicolumn{1}{|c|}{RNN}  & 9.0 s & 2,431 kg& 69.3  g& 6,482 veh

\\\hline

\multirow{2}{*}{Alternative}   & \multicolumn{1}{|c|}{FFNN} & 9.3 s & 2,432 kg& 68.7  g& 6,488 veh
\\
\multirow{1}{*}{}                & \multicolumn{1}{|c|}{RNN}      & 9.6 s & 2,430  kg& 68.6  g& 6,489 veh

\\\hline

		\end{tabular}
		\label{tabrelocation}
\end{table}

Table \ref{tabrelocation} shows the performances of the two cases of sensor location, obtained from 30 independent on-line tests. Each case is run with both FFNN and RNN as the NDR. We can see  that the two sensor locations yield very similar performances, regardless of the neural networks chosen. One possibility is that the neural networks are sufficiently deep such that their performances are not sensitive to the dimension or the nature of the state variables, as long as their parameters are sufficiently trained off line. It also shows that the proposed NDR is quite robust against different configurations of the network of sensors.

\subsubsection{Performances with demand increase}\label{subsubsecdincrease}

To further test the performance of the proposed signal controls under more congested network conditions, we increase the dynamic travel demand in the network by uniformly scaling up the demand matrix by 10\%, 20\% and 30\%. Table \ref{tabdincrease} shows the corresponding performances of FFNN and RNN (the off-line optimization minimizes traffic delay), and compares them with the baseline scenario (original demand). All figures reported are based on 30 independent simulation runs.

\begin{table}[h!]
	\centering
		\caption{On-line performances of the NDR with increased travel demand. The percentages indicate the relative increases compared to the baseline (based on the same type of neural network).}
\begin{tabular}{| c   c | l | l | l | l |}
			\hline
			 &   & Delay (s) & CO$_2$ (kg) & BC (g)  & Throughput (veh)
			 \\\hline
\multicolumn{1}{|c}{Baseline}  & \multicolumn{1}{|c|}{FFNN}  & 8.7 / - & 2,432 / -& 68.6 / - & 6,487 / - 
\\
\multicolumn{1}{|c}{(0\% increase)}  &  \multicolumn{1}{|c|}{RNN}  & 9.0 / - & 2,431 / -& 69.3 / - & 6,482 / -

\\\hline

\multirow{2}{*}{10\% increase}   & \multicolumn{1}{|c|}{FFNN} & 14.0  / 61\% & 2,627 / 8\%& 72.8 / 6\% & 7,068 / 9\%
\\
\multirow{1}{*}{}                & \multicolumn{1}{|c|}{RNN}      & 12.6 / 40\% & 2,622 / 8\% & 73.2 / 6\% & 7,062 / 9\%

\\\hline

\multirow{2}{*}{20\% increase}   & \multicolumn{1}{|c|}{FFNN} & 20.3 / 133\% & 2,860 / 18\% & 76.3 / 11\% & 7,737 / 19\%
\\
\multirow{1}{*}{}                & \multicolumn{1}{|c|}{RNN}      & 17.6 / 96\% & 2,853 / 17\%& 78.7 / 14\% & 7,743 / 19\%

\\\hline

\multirow{2}{*}{30\% increase}   & \multicolumn{1}{|c|}{FFNN} & 41.8 / 380\% & 3,219 / 32\% & 83.9 / 22\% & 8,403 / 30\%
\\
\multirow{1}{*}{}                & \multicolumn{1}{|c|}{RNN}      & 36.2 / 302\% & 3,182 / 31\%& 83.3 / 20\% & 8,448 / 30\%

\\\hline
		\end{tabular}
		\label{tabdincrease}
\end{table}

It can be seen from the table that the vehicle throughputs are consistent with the demand increase (10\%, 20\% and 30\%). Regarding the other three performance indicators, there is a hyperlinear increase of delays as the network demand increases, which is caused by the nonlinear effect of network dynamics. In comparison, CO$_2$ emissions exhibit linear growth with the demand levels, which suggests that CO$_2$ emissions are proportional to the vehicle volumes, and are less sensitive to vehicle delays. Note that this does not contradict Figure \ref{figtradeoff}, which shows high correlation between delay reduction and CO$_2$ emissions, for two reasons. Firstly, the differences of delay and CO$_2$ are due to different demand levels in Table \ref{tabdincrease}, instead of different control strategies in Figure \ref{figtradeoff}. Secondly, Figure \ref{figtradeoff} only shows local reductions at the junction level, when in fact vehicles emit majority (around 63\%) of the CO$_2$ while traveling along the links. Finally, BC emissions grow sublinearly with demand. This is quite interesting as BC emissions are mainly produced during stop-and-go cycles near junctions, which suggests that the NDR approach is effective in reducing BC emissions.

When the network becomes more congested, RNN starts to outperform FFNN in delay reduction, by 10\%, 13\% and 13\% respectively, under 10\%, 20\% and 30\% demand increase. This is because RNN takes into account the temporal precedence and chronological dependencies of the input variables when generating control parameters, and hence is capable of handling the highly nonlinear traffic dynamics under higher network loads.

\section{Conclusion}\label{secConclusion}

This paper develops a real-time signal control framework based on nonlinear decision rule (NDR) to allow actuation of signal timing changes based on network traffic states. The NDR has been implemented with two neural networks: feedforward neural network (FFNN) and recurrent neural network (RNN). Through the NDR, the controller updates traffic signal parameters based on prevailing network states, and the performance of such mechanism can be optimized via off-line training of the NDR. Particle swarm optimization is employed to solve the off-line optimization problem, which is the computationally expensive part of the NDR framework, and the on-line implementation of the trained NDR is quite efficient and can accommodate real-time decision requirements. This is a key advantage of the NDR approach.

We demonstrate the applicability and effectiveness of the proposed approaches using microscopic traffic simulation and emission modeling based on a real-world traffic network in west Glasgow. The traffic and emission models have been set up and calibrated based on an EU project (http://www.carbotraf.eu). Historical traffic flow data are used to reflect the levels of traffic demand and variability. The test phase is conducted in a simulation environment with different random seeds to populate stochasticity in the simulation. The performance of the proposed NDR approach is assessed in terms of travel delay, throughput, total carbon and black carbon emissions. The following findings are made.

\begin{itemize}
\item Compared with the fixed-timing plan used on the real-world site, the proposed NDR reduces network-wide delay by up to 68\%, total carbon and black carbon emissions by 3\% and 2\%, respectively, and 1\% increase of network throughput. In addition, most emission reductions take place at signalized intersections, as a result of the proposed controls. 

\item Under the normal network demand level, the performances of FFNN and RNN are similar in terms of delay, CO$_2$ and BC emissions, and throughput. When the network demand increases (by 10\%, 20\% and 30\% in this paper), RNN begins to outperform FFNN. This is likely due to the internal structures of FFNN and RNN as we explained at the end of Section \ref{subsubsecdincrease}. Furthermore, there seems to be a mismatch between the depth of the neural networks and the nonlinearity of the traffic/control dynamics (i.e. `depth' of the traffic network). The latter is dependent on the level of saturation of the traffic network, which causes the change in the relative performances of FFNN and RNN.

\item There is a strong correlation between delay reductions and CO$_2$ emissions at local intersections. Such a correlation does not exist between delay reductions and BC emissions. This is because CO$_2$ emissions are highly dependent on vehicle average speeds, which are related to junction delays; BC emissions, on the other hand, are affected by stop-and-go cycles and vehicle type (such as buses and HGVs), which are not directly related to junction delays. Furthermore, minimizing delays in the off-line training tends to also minimize CO$_2$ and BC emissions.

\item The NDR approaches with FFNN and RNN are both tested with a different set of loop detectors in the network, which offers relatively more complete information on all the incoming approaches of signalized intersections. Note that in the real-world network, some intersections have missing detectors on some incoming approaches. The test result, surprisingly, shows that the performances are very similar in these two cases, which means that the NDR approach is robust against different sensor locations. 
\end{itemize}

Further extensions of the proposed neural-network-based NDR approach include 
\begin{itemize}
\item[(1)] integration with localized signal coordination (such as offset optimization) and actuation (such a bus signal priority) to further reduce emissions at junctions; 
\item[(2)] a systematic and quantitative approach to optimal sensor location that is compatible with the proposed neural network structure; and 
\item[(3)] an in-depth investigation of the matching between the structure and depth of the neural networks, and the depth of the traffic network. Such a priori knowledge could offer insights into the construction of nonlinear decision rules and tuning of their parameters for better efficiency and less redundancy.
\end{itemize}

\section*{Acknowledgement}
This study was partially supported by the CARBOTRAF project (no. 28786) funded by the EU 7th Framework Program, National Social Science Foundation of China (15BGL143), and the Zhejiang University/University of Illinois at Urbana-Champaign Institute.


\begin{thebibliography}{99}


\bibitem[Abu-Mostafa, 1989]{ABUM1989} Abu-Mostafa, Y.S., 1989. Information theory, complexity and neural networks. IEEE Communications Magazine, 27(11), pp.25-28.

\bibitem[Arel et al., 2010]{ALUK2010} Arel, L., Liu, C., Urbanik, T., Kohls, A.G., 2010. Reinforcement learning-based multi-agent system for network traffic signal control, IET Intelligent Transport Systems, 4(2), 128-135.



\bibitem[Balaji et al., 2010]{BGS2010} Balaji P.G., German X., Srinivasan. D., 2010. Urban traffic signal control using reinforcement learning agents, IET Intelligent Transport Systems, 4(3), 177-188.


\bibitem[Banks et al., 2007]{BVA2007} Banks, A., Vincent, J., Anyakoha, C., 2007. A review of particle swarm optimization. Part I: background and development. Nat. Comput. 6 (4), 467-484.


\bibitem[Bertsimas et al., 2011]{BBC2011} Bertsimas, D., Brown, D.B., Caramanis, C., 2011. Theory and applications of robust optimization. SIAM Rev. 53 (3), 464-501.


\bibitem[Cai et al., 2009]{CWH2009} Cai, C., Wong, C.K., Heydecker, B.G., 2009. Adaptive traffic signal control using approximate dynamic programming. Transportation Research Part C: Emerging Technologies, 17(5), pp.456-474.


\bibitem[Castro et al., 2017]{CHM2017} Castro, G.B., Hirakawa, A.R., Martini, J.S., 2017. Adaptive traffic signal control based on bio-neural network. Procedia Computer Science, 109, 1182-1187.

\bibitem[Chang and Hui, 2016]{CH2016} Chang, L., Hui, W., 2016, Traffic emission control based on emission pricing and signal timing. In Intelligent Control and Automation (WCICA), 2016 12th World Congress on (pp. 467-472). IEEE.


\bibitem[Chang and Sun, 2004]{CS2004} Chang, T.H., Sun, G.Y., 2004. Modeling and optimization of an oversaturated signalized network. Transportation Research Part B, 38(8),  687-707.


\bibitem[Chen et al., 2012]{CBMW2012} Chen, H., Bai, R., Ma, J., Wang, D, 2012. Research on Intersection Signal Timing Model Considering Emissions Effects. In CICTP 2012, American Society of Civil Engineers, 1024-1034.


\bibitem[Christofa et al., 2016]{CAS2016} Christofa, E., Ampountolas, K., Skabardonis, A., 2016. Arterial traffic signal optimization: A person-based approach. Transportation Research Part C, 66, 27-47.


\bibitem[Christofa and Skabardonis, 2011]{CS2011} Christofa, E., Skabardonis, A., 2011. Traffic signal optimization with application of transit signal priority to an isolated intersection. Transportation Research Record: Journal of the Transportation Research Board, 2259,192-201.


\bibitem[Feng et al., 2015]{FHKZ2015} Feng, Y., Head, K.L., Khoshmagham, S., Zamanipour, M., 2015. A real-time adaptive signal control in a connected vehicle environment. Transportation Research Part C, 55, 460-473.


\bibitem[Friesz, 2010]{Friesz2010} Friesz, T.L., 2010. Dynamic Optimization and Differential Games. Springer, New York.

\bibitem[Fu and chen, 1993]{FUCHEN1993} Fu, L. and Chen, T., 1993. Sensitivity analysis for input vector in multilayer feedforward neural networks. In Neural Networks, 1993., IEEE International Conference on (pp. 215-218). IEEE.

\bibitem[Gartner, 1983]{OPAC} Gartner, N.H., 1983. OPAC: a demand-responsive strategy for traffic signal control. Transportation Research Record. 906: 75-81.

\bibitem[Gkatzoflias, et al., 2006]{GKNS2006} Gkatzoflias D, Kouridis C, Ntziachristos L, Samaras Z., 2006. COPERT 4 manual. European Environment Agency (EEA)


\bibitem[Han, 2017]{Han2017} Han, K, 2017. Framework for real-time traffic management with case studies. Transportation Research
Record: Journal of the Transportation Research Board, No. 2658, 35-43.


\bibitem[Han and Gayah, 2015]{HG2015} Han, K, Gayah, VV, 2015. Continuum signalized junction model for dynamic traffic networks: Offset, spillback, and multiple signal phases. Transportation Research Part B 77, 213-239.


\bibitem[Han et al., 2014]{HGPFY2014} Han, K., Gayah, V.V., Piccoli, B., Friesz, T.L., Yao, T., 2014. On the continuum approximation of the on-and-off signal control on dynamic traffic networks. Transportation Research Part B, 61, 73-97.


\bibitem[Han et al., 2016]{HLGFY2016} Han, K, Liu, H, Gayah, VV, Friesz, TL, Yao, T., 2016. A robust optimization approach for dynamic traffic signal control with emission considerations. Transportation Research Part C, 70, 3-26.


\bibitem[Han et al., 2015]{HSLFY2015} Han, K, Sun, Y, Liu, H, Friesz, TL, Yao, T, 2015. A bi-level model of dynamic traffic signal control with continuum approximation. Transportation Research Part C, 55, 409-431.


\bibitem[Hauser and Scherer, 2001]{HS2001} Hauser. T.A., Scherer. W.T., 2001. Data mining tools for real-time traffic signal decision support \& Maintenance, Systems, Man, and Cybernetics, 2001 IEEE International Conference on, 3, 1471-1477.


\bibitem[He et al., 2014]{HHD2014} He, Q., Head, K.L., Ding, J., 2014. Multi-modal traffic signal control with priority, signal actuation and coordination. Transportation Research Part C, 46, 65-82.


\bibitem[Henry, 1983]{PRODYN} Henry, J.J., Farges, J.L., Tuffal, J., 1983. The PRODYN real time traffic algorithm. Proceedings of the fourth IFAC-IFIP-IFORS conference on Control in Transportation Systems. 307-311.


\bibitem[Hunt et al., 1982]{SCOOT} Hunt, P.B., Robertson, D.I., Bretherton, R.D., 1982. The SCOOT on-line traffic signal optimization technique. Traffic Engineering and Control. 25, 14-22.


\bibitem[Jamshidnejad et al., 2017]{JPPD2017} Jamshidnejad, A., Papamichail, I., Papageorgiou, M. and De Schutter, B., 2017. Sustainable model-predictive control in urban traffic networks: Efficient solution based on general smoothening methods. IEEE Transactions on Control Systems Technology, 26(3), 813-827. 


\bibitem[Janssen et al., 2013]{JGLSCHFBK2013} Janssen NAH, Gerlofs-Nijland ME, Lanki T, Salonen RO, Cassee F, Hoek G, Fischer P, Brunekreef B,
Krzyzanowski M (2013) Health effects of black carbon. In: Europe WROf (ed)

\bibitem[Ji et al., 2014]{JHHT2014} Ji, Y., Hu, B., Han, J., Tang, D., 2014. An improved algebraic method for transit Signal priority scheme and its impact on traffic emission. Mathematical Problems in Engineering, 2014, 11 pages.

\bibitem[Lefebvre et al., 2011]{LEFE2011} Lefebvre, W., F. Fierens, E. Trimpeneers, S. Janssen, K. Van de Vel, F. Deutsch, P. Viaene, J. Vankerkom, G. Dumont, C. Vanpoucke, C. Mensink, W. Peelaerts, and J. Vliegen. Modeling the effects of a speed limit reduction on traffic-related elemental carbon (EC) concentrations and population exposure to EC. Atmospheric Environment, Vol. 45, No. 1, 2011, pp. 197-207.

\bibitem[Li et al., 2016]{LLW2016} Li, L., Lv, Y., Wang, F., 2016. Traffic Signal Timing via Deep Reinforcement Learning, IEEE/CAA Journal or Automatica Sinica, 3(3), 247-254.


\bibitem[Lin et al., 2013]{LDXH2013} Lin, S., De Schutter, B., Xi, Y., Hellendoorn, H., 2013. Integrated urban traffic control for the reduction of travel delays and emissions, IEEE Transactions on Intelligent Transportation Systems, 14(4), 1609-1619.


\bibitem[Liu et al., 2015]{LHGFY2015} Liu, H., Han, K., Gayah, V.V., Friesz, T.L., Yao, T., 2015. Data-driven linear decision rule approach for distributionally robust optimization of on-line signal control. Transportation Research Part C, 59, 260-277.


\bibitem[Lowrie, 1982]{SCAT} Lowrie, P., 1982. The Sydney coordinated adaptive traffic system-principles, methodology, algorithms. International Conference on Road Traffic Signalling.


\bibitem[Lucas et al., 2000]{LMH2000} Lucas, D., Mirchandani, P., Head, K., 2000. Remote simulation to evaluate real-time traffic control strategies. Transportation Research Record: Journal of the Transportation Research Board, 1727, 95-100.



\bibitem[Mascia et al., 2015]{MHHN2015} Mascia, M., Hu, K., Han, K., North, R.,  2015. Simulation output for traffic scenarios for the city of Glasgow. Technical report, CARBOTRAF, D3.4.  




\bibitem[Mascia et al., 2016]{MHHNVTBL2016} Mascia, M., Hu, K., Han, K., North, R., Van Poppel, M., Theunis, J., Beckx, C.,  Litzenberger, M., 2016. Impact of Traffic Management on Black Carbon Emissions: a Microsimulation Study. Networks and Spatial Economics, 17(1), 269-291.


\bibitem[Osorio et al., 2015]{OSO2015} Osorio, C. and Nanduri, K., 2015. Urban transportation emissions mitigation: Coupling high-resolution vehicular emissions and traffic models for traffic signal optimization. Transportation Research Part B: Methodological, 81, pp.520-538.


\bibitem[Papatzikou and Stathopoulos, 2015]{PS2015} Papatzikou, E., Stathopoulos, A., 2015. An optimization method for sustainable traffic control in urban areas. Transportation Research Part C, 55, 179-190.


\bibitem[Samah et al., 2013]{SBH2013} Samah, E., Baher, A., Hossam A., 2013. Multiagent reinforcement learning for integrated network of adaptive traffic signal controllers (MARLIN-ATSC): Methodology and Large-scale application on downtown Toronto, IEEE transactions on Intelligent transportation systems, 14(3), 1140-1150.

\bibitem[Savsani et al., 2010]{SRV2010} Savsani, V., Rao, R.V. and Vakharia, D.P., 2010. Optimal weight design of a gear train using particle swarm optimization and simulated annealing algorithms. Mechanism and machine theory, 45(3), 531-541.

\bibitem[Sha and Hsu, 2008]{SH2008} Sha, D.Y. and Hsu, C.Y., 2008. A new particle swarm optimization for the open shop scheduling problem. Computers \& Operations Research, 35(10), pp.3243-3261.

\bibitem[Smith et al., 2013]{SBXR2013} Smith, S., Barlow, G., Xie, X.-F., Rubinstein, Z., 2013. SURTRAC: Scalable Urban Traffic Control. Transportation Research Board 92nd Annual Meeting, Jan. 2013. 


\bibitem[Sobrino et al., 2016]{SMH2016} Sobrino, N., Monzon, A., Hernandez, S.,  2016. Reduced carbon and energy footprint in highway operations: the
Highway Energy Assessment (HERA) methodology. Networks and Spatial Economics, 16, 395-414.


\bibitem[S-Paramics, 2011]{Paramics2011} S-Paramics. S-Paramics 2011 reference manual.In, SIAS limited, Edinburg, Scotland, 2011.


\bibitem[Srinivasan et al., 2006]{SCC2006} Srinivasan, D., Choy, M.C., Cheu, R.L., 2006. Neural networks for real-time traffic signal control. IEEE Trans-actions on Intelligent Transportation Systems, 7(3), pp.261-272.



\bibitem[Stevanovic et al., 2015]{SSSO2015} Stevanovic, A., Stevanovic, J., So, J., Ostojic, M, 2015. Multi-criteria optimization of traffic signals: Mobility, safety, and environment. Transportation Research Part C, 55, 46-68.



\bibitem[Sun et al., 2006]{SBW2006} Sun, D., Benekohal, R.F., Waller, S.T., 2006. Bi-level programming formulation and heuristic solution approach for dynamic traffic signal optimization. Computer-Aided Civil and Infrastructure Engineering, 21(5), 321-333.


\bibitem[Sundaram et al., 2015]{SKD2015} Sundaram, S., Kumar, S.S., Divya Shree, M.S., 2015. Hierarchical clustering technique for traffic signal decision support. International Journal of Innovative Science, 2(6), 72 - 82.


\bibitem[Sunkari, 2004]{Sunkari} Sunkari, S., 2004. The benefits of retiming traffic signals, Institute of Transportation Engineers, ITE Journal, 74(4) 26 pages. 


\bibitem[AIRE, 2011]{AIRE2011} Transport ScotlandÕs Instantaneous Emissions Software AIRE. 2011. http://www.sias.com/2013/AIRE.htm.

\bibitem[Transport Scotland, 2011]{TS2011}. Transport Scotland. AIRE (Analysis of Instantaneous Road Emissions) User Guidance In, Scotland, 2011.

\bibitem[Ukkusuri et al., 2010]{URP2010} Ukkusuri, S.V., Ramadurai, G., Patil, G., 2010. A robust transportation signal control problem accounting for traffic dynamics. Computers and Operations Research 37 (5), 869-879.


\bibitem[U.S. Environmental Protection Agency, 2012]{USEPA2012} U.S. Environmental Protection Agency (2012) Report to Congress on Black Carbon. Department of the Interior, Environment, and Related Agencies Appropriations Act. EPA-450/R-12-001

\bibitem[Webster, 1958]{Webster1958} Webster, F. V. Traffic Signal Settings. H.M. Stationery Office, 1958.

\bibitem[Wiering, 2000]{Wiering2000} Wiering, M.A., 2000. Multi-agent reinforcement learning for traffic light control. In Proceedings of the 17th International Conference on Machine Learning, 1151-1158.

\bibitem[Yin, 2006]{Yin2006} Yin, P.Y., 2006. Particle swarm optimization for point pattern matching. Journal of Visual Communication and Image Representation, 17(1),143-162.

 \bibitem[Yin, 2008]{Yin2008} Yin, Y., 2008. Robust optimal traffic signal timing. Transportation Research Part B, 42(10), 911-924.
 
 
 \bibitem[Zhang et al., 2011]{ZBD2011} Zhang, K., Batterman, S., Dion, F., 2011. Vehicle emissions in congestion: Comparison of work zone, rush hour and free-flow conditions. Atmospheric Environment, 45(11), 1929-1939.
 
 
\bibitem[Zhang et al., 2010]{ZYL2010} Zhang, L., Yin, Y., Lou, Y., 2010. Robust signal timing for arterials under day-to-day demand variations. Transportation Research Record: Journal of the Transportation Research Board, 2192, 156-166.
 
 
 \bibitem[Zhang et al., 2013]{ZYC2013} Zhang, L., Yin, Y., Chen, S., 2013. Robust signal timing optimization with environmental concerns. Transportation Research Part C 29, 55-71.
 
 \bibitem[Zhou and Cai, 2014]{ZC2014} Zhou, Z., Cai, M., 2014. Intersection signal control multi-objective optimization based on genetic algorithm. Journal of Traffic and Transportation Engineering (English Edition), 1(2), 153-158.
 

\end{thebibliography}
\end{document}